\documentclass[final, 3p, 11pt, authoryear]{elsarticle}

\usepackage{amsmath,amssymb,amsfonts}

\usepackage{amsmath,amssymb,amsfonts}
\usepackage[dvipsnames]{xcolor}
\definecolor{refcolor}{RGB}{40, 56, 104}
\usepackage{hyperref}
\hypersetup{
  breaklinks=true,
  colorlinks=true,
  citecolor=blue,
  linkcolor=blue,
}
\usepackage[ruled, vlined, boxed]{algorithm2e}
\usepackage{booktabs}
\usepackage{multirow}
\usepackage{makecell}
\usepackage{threeparttable}

\journal{arXiv}

\begin{document}

\begin{frontmatter}



\title{Deep Reinforcement Learning for Traveling Purchaser Problems}

\author[1]{Haofeng Yuan}
\ead{yhf22@mails.tsinghua.edu.cn}

\author[1]{Rongping Zhu}
\ead{zhurp19@mails.tsinghua.edu.cn}

\author[1]{Wanlu Yang}
\ead{yangwl21@mails.tsinghua.edu.cn}

\author[1]{Shiji Song\corref{cor1}}
\ead{shijis@mail.tsinghua.edu.cn}

\author[1]{Keyou You}
\ead{youky@tsinghua.edu.cn}

\author[2]{Wei Fan}
\ead{fanwei2@lenovo.com}

\author[3]{C.~L.~Philip~Chen}
\ead{philip.chen@ieee.org}

\cortext[cor1]{Corresponding author}

\address[1]{Department of Automation \& BNRist, Tsinghua University, Beijing 100084, China}
\address[2]{AI Laboratory, Lenovo Research, Beijing 100193, China}
\address[3]{School of Computer Science and Engineering, South China University of Technology, Guangzhou 510641, China}

\begin{abstract}
The traveling purchaser problem (TPP) is an important combinatorial optimization problem with broad applications. 
Due to the coupling between routing and purchasing, existing works on TPPs commonly address route construction and purchase planning simultaneously, which, however, leads to exact methods with high computational cost and heuristics with sophisticated design but limited performance.
In sharp contrast, we propose a novel approach based on deep reinforcement learning (DRL), which addresses route construction and purchase planning \textit{separately}, while evaluating and optimizing the solution from a \textit{global} perspective.
The key components of our approach include a bipartite graph representation for TPPs to capture the market-product relations, and a policy network that extracts information from the bipartite graph and uses it to sequentially construct the route.
\textcolor{black}{One signif\mbox{}icant advantage of our framework is that we can eff\mbox{}iciently construct the route using the policy network, and once the route is determined, the associated purchasing plan can be easily derived through linear programming, while, by leveraging DRL, we can train the policy network towards optimizing the global solution objective.}
Furthermore, by introducing a meta-learning strategy, the policy network can be trained stably on large-sized TPP instances, and generalize well across instances of varying sizes and distributions, even to much larger instances that are never seen during training.
Experiments on various synthetic TPP instances and the TPPLIB benchmark demonstrate that our DRL-based approach can signif\mbox{}icantly outperform well-established TPP heuristics, reducing the optimality gap by 40\%-90\%, and also showing an advantage in runtime, especially on large-sized instances.

\end{abstract}


\begin{keyword} 
Traveling purchaser problem \sep deep reinforcement learning \sep bipartite graph \sep policy network \sep meta-learning
\end{keyword}

\end{frontmatter}

\vfill\newpage
\setcounter{tocdepth}{2}
\tableofcontents
\vfill\newpage

\section{Introduction}

The traveling purchaser problem (TPP) is a well-known combinatorial optimization (CO) problem, which has broad real-world applications and received much attention from both researchers and practitioners in recent decades \citep{r_49, r_21, r_43, r_44}. In the TPP, given a list of products and their demand quantities, the purchaser aims to decide a route and an associated purchasing plan to meet the product demands by visiting a subset of markets, with the objective to minimize the sum of traveling and purchasing costs \citep{r_23}.

Unfortunately, TPPs are known to be strongly NP-hard \citep{r_1}. In particular, the main challenge in solving TPPs stems from the inherent coupling between routing and purchasing: the routing decisions have to take both traveling and purchasing costs into consideration, while the choice of visited markets will, in turn, impact the purchasing decisions and associated costs. Thus, existing works typically need to deal with both routing and purchasing simultaneously \citep{r_1}, which, however, imposes many limitations. On the one hand, exact methods such as branch-and-cut \citep{r_22, r_41} require joint optimization of routing and purchasing decisions, which generically leads to high computational cost, making them often intractable for realistically sized problem instances. For example, solving a medium-sized TPP instance with 100 markets and 50 products can take hours on a computer with a 2.3 GHz processor, and the computational time grows exponentially with the problem size \citep{r_35}. On the other hand, various heuristic methods, such as the Generalized Savings Heuristic (GSH) \citep{r_24}, the Tour Reduction Heuristic (TRH) \citep{r_25}, and the Commodity Adding Heuristic (CAH) \citep{r_26}, have been proposed to produce high-quality solutions within a reasonable time. 
Readers are referred to the Appendices for further details.
\textcolor{black}{However, these heuristics have to carefully balance the effect on traveling and purchasing costs for each operation, and thus rely heavily on sophisticated design that requires substantial human expertise and domain knowledge.
Moreover, these hand-crafted heuristics typically need to be tailored case-by-case, resulting in that they may only be effective on limited instances with specific characteristics, and lack the ability to generalize across different instance distributions.} For example, CAH performs poorly if the purchasing cost dominates in the optimization objective, and GSH may yield poor solutions if a market sells most products but is located far from other markets \citep{r_26}. 

In this paper, we propose a novel approach based on deep reinforcement learning (DRL), which addresses the limitations of existing methods by decoupling the treatment of routing and purchasing decisions. The core idea of our approach can be summarized as ``\emph{solve separately, learn globally}''.
We notice that by leveraging the forward-thinking and global-optimizing mechanisms of DRL, we can break the complex task of solving TPPs into two \textit{separate} stages: route construction and purchase planning, while learning a policy to guide the decision-making towards optimizing the \textit{global} solution objective. 
Specifically, in the first stage, we use a policy network to sequentially construct the route, where at each decision step the policy network takes the problem instance and current partial route as input, and outputs an action distribution to determine the next market to visit. %
Then, once a complete route is constructed, we proceed to the second stage, where we derive the optimal purchasing plan for the visited markets through an easily solvable linear transportation problem. %
Note that in the first stage, we use the policy network only to decide the route, i.e., which markets to visit and the visited order, deferring the decisions on specific purchasing plan to the second stage.
\textcolor{black}{This separation decouples the routing decisions and purchasing decisions at the operational level, such that each stage of decisions can be eff\mbox{}iciently done. 
Meanwhile, to bridge the optimization interdependence between routing and purchasing, we train the first-stage policy network to optimize the global solution objective, which is determined jointly by both stages.} 
In other words, though the policy network is only used for route construction in the first stage, it is learned to construct a high-quality route that can not only have low traveling cost, but also be promising to lead to a low-cost purchasing plan in the second stage, thus to minimize the global solution objective.

In fact, the idea of leveraging DRL for CO problems is not entirely new though \citep{r_28, r_54, r_59, r_55}, particularly in routing problems such as the traveling salesman problem (TSP) \citep{r_3, r_32, r_37, r_38} and the vehicle routing problem (VRP) \citep{r_2, r_46, r_48, r_76}. However, these efforts typically rely on problem-specific properties and are limited to problems with simple structures, which cannot be readily extended to TPPs. \textcolor{black}{Therefore, considering the double nature of TPPs, we first propose a novel bipartite graph representation as an input to the policy network}, where the markets and products are represented as two sets of nodes, with the supply information (e.g., the supply quantities and prices) encoded as edge features between the market and product nodes. In contrast to existing representations for routing problems (e.g., the complete or \emph{k-NN} graphs \citep{r_4, r_5}), our bipartite graph representation can naturally capture the relations between markets and products through message-passing along edges.  
Then, we design a policy network, with an architecture that can effectively extract information from the bipartite graph and use it for route construction. Specifically, our policy network exploits the connecting structure of the bipartite graph, and leverages graph neural networks (GNNs) and multi-head attention (MHA) to aggregate information both within and across the two sets of nodes. This effectively facilitates the extraction of relational information of the markets and products, such as the spatial relations between markets and the potential substitutions and complementarities in product supply, which is important for the construction of high-quality routes. \textcolor{black}{Moreover, it is worth mentioning that our bipartite graph representation and policy network are size-agnostic, such that a trained policy network can be flexibly adapted to TPP instances of varying sizes without the need for retraining or parameter adjustment.}

In addition, despite our proposed framework and the key components introduced above, how to efficiently learn an effective policy remains another challenge. Due to the huge state-action space, a randomly initialized policy network may suffer from inefficient exploration, which can result in slow convergence or even training collapse, particularly on large-sized instances. 
Different from the training techniques or tricks designed for specific problems (mainly TSP or VRP) \citep{r_36, r_45} , we introduce an effective and general training strategy based on meta-learning \citep{r_18}, which trains an initialized policy network on a collection of varying instance distributions, with the learning objective to achieve efficient adaptation to new instances at low fine-tuning cost. Furthermore, the initialized policy network can learn cross-distribution knowledge through meta-learning, such that it can effectively generalize across varied instance sizes and distributions, even demonstrating zero-shot generalization ability to much larger instances that are never seen during training.

We empirically evaluate our DRL-based approach on two types of TPPs: the restricted TPP (R-TPP), where the supply quantities of available products at each market are limited, and the unrestricted TPP (U-TPP), which assumes unlimited supply quantities. We conduct extensive experiments on various synthetic TPP instances and the TPPLIB benchmark \citep{r_35}. The results demonstrate that our DRL-based approach can produce near-optimal solutions with high computational efficiency, significantly outperforming well-established TPP heuristics in both solution quality and runtime. Notably, our approach achieves an average optimality gap consistently within 6\% on each category of the TPPLIB benchmark in our experiment, yielding a reduction of 40\%-90\% compared to the baseline heuristics. In addition, we further confirm the zero-shot generalization ability of the learned policy on instances that are much larger (up to 300 markets and 300 products, the largest size in TPPLIB) than the training instances.

The remainder of this paper is organized as follows. In Section \ref{Problem Formulation}, we formally describe TPPs and introduce the mathematical formulation. Section \ref{MDP Reformulation} presents our ``solve separately, learn globally'' framework for solving TPPs. Section \ref{Bipartite Graph Representation for TPPs} introduces the bipartite graph representation for TPPs, and Section \ref{Policy Network} presents the details of our policy network, followed by a description of the basic training algorithm. In Section \ref{Meta-learning Stategy}, we propose a meta-learning strategy for eff\mbox{}iciently training on large-sized problems and improving generalization. Section \ref{Experiments} provides the empirical evaluation and reports the results and analysis. Finally, conclusions are drawn in Section \ref{Conclusion}.

\section{Problem Formulation} \label{Problem Formulation}

In the TPP, a purchaser needs to buy a set of products from a set of markets. Each product can only be purchased from certain markets, with potentially varying supply quantities and prices at different markets. The purchaser aims to decide a route that visits a subset of markets and an associated purchasing plan to meet all product demands, with the objective to minimize the sum of traveling and purchasing costs.

Mathematically, a TPP is defined on a complete directed graph $\mathcal{G} = \left(V, E\right)$, where $V = \left\{v_0\right\} \cup M$ is the set of nodes and $E = \{(i, j): i, j \in V, i \neq j\}$ is the set of edges. Node $v_0$ denotes the depot and $M$ denotes the set of markets. The problem considers a set $K$ of products to purchase, where a demand $d_k$ is specified for each product $k \in K$. Each product $k$ is supplied in a subset of markets $M_k \subseteq M$, where at most $q_{ik}$ units of product $k$ can be purchased from market $i \in M_k$, at a price of $p_{ik}$. The traveling cost associated with each edge $(i, j) \in E$ is denoted by $c_{ij}$. The goal is to determine a route on $\mathcal{G}$, i.e., a simple cycle $\pi = \left(v_0, v^1, v^2, \dots, v^{T}, v_0 \right),\; T \le |M|$, which starts and ends at the depot and visits a subset of markets, and decide how much of each product to purchase at each market to satisfy the demands at minimum traveling and purchasing costs. Note that to guarantee the existence of a feasible purchasing plan, it is assumed that $0 < q_{ik} \le d_k$ for each $k \in K$ and $i \in M_k$, and $\sum_{i \in M_k}q_{ik} \ge d_k$ for each $k \in K$. In the case of U-TPP, the first assumption becomes $q_{ik} = d_k$ for all $k \in K$ and $i \in M_k$, i.e., the supply quantities for the available products at each market are unlimited. Otherwise, the problem is referred to as R-TPP.

The TPP can be formulated as a mixed integer linear programming (MILP) problem. Let $y_i$ be a binary variable taking value $1$ if market $i$ is selected, and $0$ otherwise. Let $x_{ij}$ be a binary variable taking value $1$ if edge $(i, j)$ is traversed, and $0$ otherwise. Let $z_{ik}$ be a variable representing the quantity of product $k$ purchased at market $i \in M_k$. For a subset of nodes $V^\prime \subset V$, we define:
\begin{align*}
\delta^{+}\left(V^{\prime}\right):=\left\{(i, j) \in E: i \in V^{\prime}, j \notin V^{\prime}\right\}, \\
\delta^{-}\left(V^{\prime}\right):=\left\{(i, j) \in E: i \notin V^{\prime}, j \in V^{\prime}\right\},
\end{align*}
\textcolor{black}{where $\delta^{+}\left(V^{\prime}\right)$ and $\delta^{-}\left(V^{\prime}\right)$ denote the edges going out of and into the given node set $V^{\prime}$, respectively.} Then, the TPP can be formulated as follows \citep{r_22}:
\begin{align}
\text{min} & \quad \hspace{-0.3em} \sum_{(i, j) \in E} c_{ij}x_{ij} + \sum_{k \in K} \sum_{i \in M_k} p_{ik}z_{ik} \label{eq_1}\\
\text{s.t.} \hspace{0.075em} & \quad \sum_{i \in M_k} z_{ik} = d_k, \; \forall k \in K,  \label{eq_2}\\
 & \quad z_{ik} \le q_{ik}y_i, \; \forall k \in K, i \in M_k,  \label{eq_3}\\
 & \quad \hspace{-1.275em} \sum_{(i, j) \in \delta^{+}(\{h\})} x_{ij} = y_{h}, \; \forall h \in V, \label{eq_4}\\
 & \quad \hspace{-1.275em} \sum_{(i, j) \in \delta^{-}(\{h\})} x_{ij} = y_{h}, \; \forall h \in V,  \label{eq_5}\\
  & \quad \hspace{-1.275em} \sum_{(i, j) \in \delta^{-}\left(M^{\prime}\right)} x_{i j} \geq y_{h}, \; \forall M^{\prime} \subseteq M, h \in M^{\prime},  \label{eq_6}\\
 & \quad x_{ij} \in \left\{0, 1\right\}, \; \forall (i, j) \in E,  \label{eq_7}\\
 & \quad y_{i} \in \left\{0, 1\right\}, \; \forall i \in V,  \label{eq_8}\\
 & \quad z_{ik} \ge 0, \; \forall k \in K, i \in M_k. \label{eq_9}
\end{align}
The objective function (\ref{eq_1}) aims at minimizing the sum of traveling and purchasing costs. Constraints (\ref{eq_2}) ensure that all demands must be satisfied. Constraints (\ref{eq_3}) limit the quantity of product $k$ purchased at market $i$ to the available supply quantity, and prevent any purchase at unvisited markets. Constraints (\ref{eq_4})-(\ref{eq_5}) impose that for each visited market, exactly one edge must enter and leave the node, and constraints (\ref{eq_6}) prevent sub-tours. Constraints (\ref{eq_7})-(\ref{eq_9}) impose integrality conditions and bounds on the decision variables.

The MILP formulation (\ref{eq_1})-(\ref{eq_9}) is commonly employed in exact methods for TPPs \citep{r_1}. However, the number of constraints is exponential in the number of markets, and the number of variables is also significantly large, which can lead to a massive branch-and-bound tree with weak lower bounds in an MILP solver. The branch-and-cut algorithms, as introduced in \cite{r_22} and \cite{r_41}, attempt to alleviate this by dynamically generating variables and separating constraints, thereby to reduce the tree size and accelerate solving. Nonetheless, the dynamic pricing and separation procedures are still computationally expensive, making it often unaffordable to exactly solve the problem in practice, especially for large-sized TPP instances. In contrast, our DRL-based approach does not directly solve the MILP, but instead creates a bipartite graph representation based on its structure (Section \ref{Bipartite Graph Representation for TPPs}), which is further used as input to the policy network for route construction.

\section{``Solve Separately, Learn Globally'' Framework}  \label{MDP Reformulation}

We introduce our ``solve separately, learn globally'' framework in this section. As aforementioned, prior works typically require addressing routing and purchasing simultaneously \citep{r_1}. In contrast, by leveraging DRL, we can decouple the routing decisions and purchasing decisions at the operational level, i.e., ``solve separately'', while guiding the decision-making using a policy that is learned to optimize the global solution objective, i.e., ``learn globally''. This framework offers the potential for both high computational efficiency and solution quality, provided the policy is well-trained.

\begin{figure*}[t]
\centering
\includegraphics[width=\textwidth, keepaspectratio]{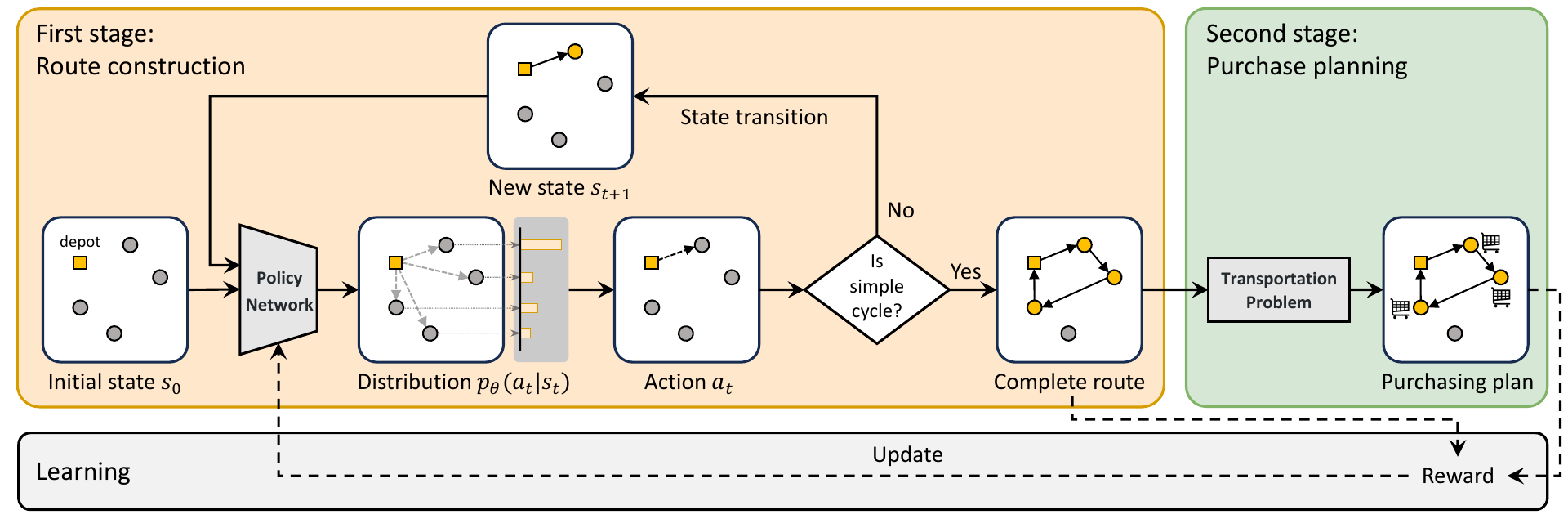}\\
\caption{Our ``solve separately, learn globally'' framework. In the first stage, the purchaser starts from the depot, and selects a new node to add to the partial route at each decision step until it returns to the depot. Then in the second stage, the purchasing plan is derived using a transportation problem, which is integrated into the reward calculation procedure for the terminal state. The final reward, defined as the negative solution objective value, is used to update the policy network through a DRL training algorithm.}
\label{fig_1}
\end{figure*}

\subsection{Overview}

Figure \ref{fig_1} illustrates an overview of our proposed framework. The solving procedure consists of two separate stages: \emph{route construction} and \emph{purchase planning}. \textcolor{black}{In implementation, we formulate this process as a Markov decision process (MDP), based on which we use DRL to guide the decisions.} Specifically, in the first stage, we sequentially construct the route through finite decision steps. The purchaser starts from the depot as the initial state, and at each decision step selects the next market (or the depot) to visit based on a policy, which is parameterized as a deep neural network (called the policy network). The route construction process iterates until the depot is revisited, indicating a complete route is formed, and the MDP reaches the terminal state. Then, we proceed to the second stage, where we derive the optimal purchasing plan for the visited markets by solving a linear transportation problem. Within the MDP, the second stage is integrated into the reward calculation procedure for the terminal state. The final reward is defined as the negative solution objective value, i.e., the sum of the traveling cost of the route constructed in the first stage and the purchasing cost derived in the second stage. The reward is then used to update the first-stage policy network through a DRL training algorithm.

\subsection{MDP Formulation}

In this subsection, we introduce the details of our framework in the form of MDP.

\subsubsection{State} \label{State}

A state should contain necessary information that the policy network needs for decision-making at the current step. In our MDP, the state at decision step $t$, denoted as $s_t$, consists of two parts: 1) the TPP instance $U$ being solved, which is \emph{static} throughout the solving process, and 2) the \emph{dynamic} information of the partial route that has been constructed up to step $t$. Specifically, the TPP instance $U$ is represented as a bipartite graph, which will be introduced in detail in Section \ref{Bipartite Graph Representation for TPPs}. For the dynamic information, we consider two contents that are critical to the decision task at step $t$:
\begin{itemize}
    \item the partial route $\pi^t = \left(v_0, v^1, v^2, \dots, v^{t-1} \right)$,
    \item the remaining demand $d^t = \left(d_1^t, d_2^t, \dots, d_{\left|K\right|}^t\right)$,
\end{itemize}
where $v^1, v^2, \dots, v^{t-1}$ denote the nodes selected at previous decision steps $1, 2, \dots, t-1$, respectively, and $d_k^t$ denotes the remaining demand of product $k \in K$, assuming all possible quantities of product $k$ supplied in the partial route have been purchased. 
In the initial state $s_1$, the route starts from the depot $v_0$ and no products have been purchased, i.e., where $\pi^1 = \left(v_0\right)$ and $d^1 = \left(d_1, d_2, \dots, d_{\left|K\right|}\right)$.

\textcolor{black}{We remark that, ideally, a well-trained policy network can automatically extract information about the remaining demand $d^t$ given the instance $U$ and the current partial route $\pi^t$, but we find that explicitly providing $d^t$ to the policy network can effectively reduce the computational load for decision-making and improve the solution quality. Therefore, we include the remaining demand $d^t$ as part of the state $s_t$.}

\subsubsection{Action} \label{Action}

At each decision step $t$, an action $a_t \in A_t$ is defined as selecting the next node (i.e., a new market or the depot) to visit. If a new market is selected, the purchaser will visit it and start the next step, $t+1$. Otherwise, if the depot is selected to visit again, it means that a complete route is formed, and the route construction terminates.

Note that an arbitrarily selected action may lead to infeasible solutions, either because 1) the constructed route may not be a simple cycle, or 2) the markets on the route cannot support a feasible purchasing plan that fulfills the product demands in the second stage. Therefore, to avoid infeasibility, we introduce two masking rules over the action set $A_t$ at each decision step $t$. First, to ensure a simple cycle, markets already in the partial route $\pi^{t}$ are excluded from $A_t$. In fact, it makes no sense to revisit a market from the perspective of optimization objective, since the traveling cost will increase strictly, but no savings on the purchasing cost can be made. Second, to guarantee the existence of a feasible purchasing plan, we mask the depot from $A_t$ until the remaining demand $d_k^t = 0$ for all $k \in K$. In other words, the purchaser is forbidden to return to the depot unless all product demands can possibly be satisfied from the current partial route. \textcolor{black}{These two masking rules ensure that the MDP always produces a feasible route and purchasing plan regardless of the policy.}

\subsubsection{Transition} \label{Transition}

Once an action $a_t$ is determined, the state $s_t$ will transit deterministically to a new state $s_{t+1}$. Specifically, if the next node $v^t$ is a new market, it will be added to the end of the partial route $\pi^t$:
$$\pi^{t+1} = \left(v_0, v^1, \dots, v^{t-1}, v^{t} \right).$$
Meanwhile, for each product $k \in K$, the previously unsatisfied demand $d_k^t$ will be replenished if product $k$ is available at the newly visited market $v^{t}$:
$$d_k^{t+1} = \begin{cases}
\max \left\{0, d_k^{t} - q_{v^{t}k}\right\}, \;\; \text{if } v^{t} \in M_k, \\
d_k^{t}, \;\; \text{if } v^{t} \notin M_k, \\
\end{cases}$$
where $q_{v^{t}k}$ is the available quantity of product $k$ supplied at market $v^{t}$. Otherwise, if the next node $v^t$ is the depot, a complete route is formed and the MDP terminates and then the reward calculation procedure for the terminal state is executed. Note that the TPP instance $U$, i.e., the static part of the states, remains fixed throughout the state transitions in an MDP.

\subsubsection{Reward} \label{Reward}

After the sequential route construction process as described above, the second stage---purchase planning---is executed, which is integrated into the reward calculation procedure for the terminal state. Given the constructed route, a linear transportation problem is solved to derive the optimal purchasing plan for the visited markets. In this transportation problem, each source point corresponds to a market on the route, and each destination point corresponds to a product to be purchased. The unit transportation cost from source point $i$ to destination point $k$ is exactly the cost to purchase a unit of product $k$ from market $i$ in the TPP. This transportation problem can be efficiently solved using an off-the-shelf LP solver or a polynomial-time dynamic programming algorithm \citep{r_50}. So far, a complete solution, including a route and the associated purchasing plan, has been obtained.

The objective of TPPs is to minimize the sum of traveling and purchasing costs while meeting all problem constraints. Since the feasibility of the route and the purchasing plan have been guaranteed by the masking rules (Section \ref{Action}), we define the reward for the terminal state as the negative objective value, i.e., the sum of traveling and purchasing costs, and assign a zero-reward for all intermediate steps. The reward is used to update the policy network using a DRL training algorithm, which will be described in detail in Section \ref{Training Algorithm}. \textcolor{black}{This way, the policy network is trained to optimize the global objective, that is, to construct a high-quality route that minimizes not only the traveling cost but also the subsequent purchasing cost based on it, thus to minimize the total cost.}

\subsubsection{Policy} \label{Policy}

Given the state $s_t$ at each step $t$, the action $a_t$ is decided based on a stochastic policy $p_{\theta}(a_t \mid s_t)$, which can be viewed as a distribution over the action set $A_t$. Since the transition from a state $s_t$ to the next state $s_{t+1}$ is deterministic given an action $a_t$ (Section \ref{Transition}), the joint probability of producing a route $\pi$ based on the stochastic policy $p_{\theta}(a_t \mid s_t)$ can be factorized according to the chain rule as:
$$p_{\theta}(\pi \mid U)=\prod_{t=1}^{T} p_{\theta}(a_t \mid s_t),$$
where the route is constructed in $T$ steps, and the action $a_t$ is sampled based on the policy $p_{\theta}(a_t \mid s_t)$ at each step. In our DRL-based approach, the policy is parameterized as a deep neural network $\theta$, where the input of the policy network is the state $s_t$ and the output is the policy, i.e., the action distribution $p_{\theta}(a_t \mid s_t)$. \textcolor{black}{The DRL training aims to learn a policy that can produce high-reward solutions with high probabilities and low-reward solutions with low probabilities.}

\section{Bipartite Graph Representation for TPPs} \label{Bipartite Graph Representation for TPPs}

\begin{figure*}[t]
\centering
\includegraphics[width=\textwidth, keepaspectratio]{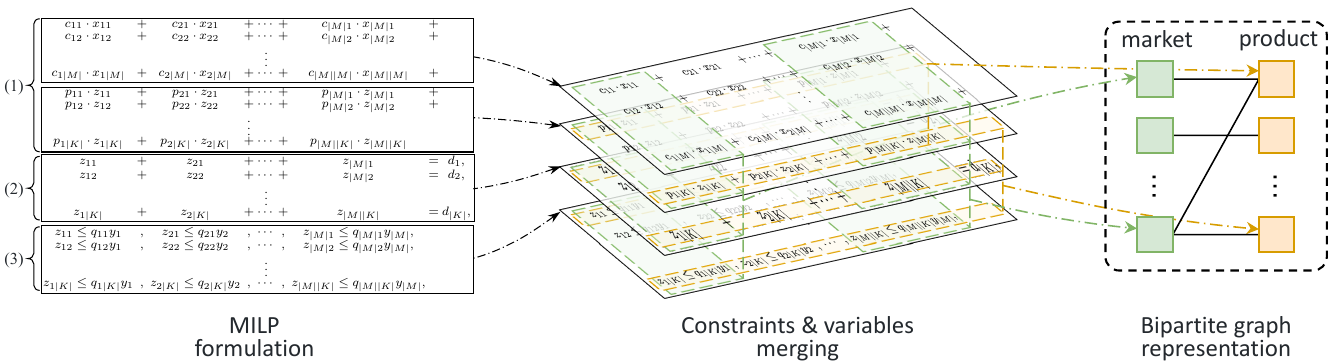}\\
\caption{The bipartite graph representation for TPPs. The bipartite graph representation is built on the MILP formulation of TPPs, where the constraints or variables that are associated with the same element are merged (as the central subfigure).}
\label{fig_4}
\end{figure*}

In the following sections, we will introduce the key DRL components in our framework. In this section, we first propose a novel bipartite graph representation for TPPs. The bipartite graph representation serves as an important part of the states (Section \ref{State}), which provides the policy network with global information about the TPP instance being solved.

Similar to existing DRL-based methods for routing problems \citep{r_4, r_5}, a natural and straightforward idea is to represent a TPP instance as a complete or \emph{k-NN} graph, where each market is represented as a node, connected to each other, and the product supply information at each market is encoded as node features. However, despite containing necessary information for the TPP instance, such representations impose many drawbacks when used in DRL. First, the feature dimensions vary with the number of products in the TPP instance being solved. Specifically, two features are needed for each product---the available quantity and the price. Consequently, a 50-product instance requires a 100-dimensional feature vector for each node to record the supply information for these 50 products, while an instance with 100 products would require 200 dimensions.
This will lead to a computational mismatch in the forward-propagation of neural networks. In other words, a policy network can only be applied on fixed-sized TPP instances, unable to adapt to instances with different sizes even if only one additional product is introduced, significantly limiting its practical usability.
Second, the relations between markets and products, such as the substitutions and complementarities in product supply between dif\mbox{}ferent markets, are very important for constructing a low-cost route. However, such dependencies are largely ignored in the complete or \emph{k-NN} graph representations. To address these limitations, we design a novel bipartite graph representation, which is built upon the MILP formulation of TPPs.

Bipartite graphs have been demonstrated effective in representing general MILP problems \citep{r_7, r_51}, where the variables and constraints are represented as two sets of nodes, and an edge connects a variable node and a constraint node if the corresponding column contributes to this constraint. Readers are referred to \cite{r_6} for a more detailed introduction. However, directly applying this representation on the MILP formulation (\ref{eq_1})-(\ref{eq_9}) for TPPs can introduce unnecessary redundancy. Specifically, the parameters to characterize a specific TPP instance appear only in the objective function (\ref{eq_1}) and constraints (\ref{eq_2})-(\ref{eq_3}), whereas constraints (\ref{eq_4})-(\ref{eq_6}), which enforce the formation of a simple cycle, are already guaranteed by the masking rules (Section \ref{Action}) and thus need not be explicitly represented for the policy network. Moreover, the constraints and variables associated with the same element (e.g., a market or a product) can be merged to reduce the graph size and thus improve the computational efficiency of the policy network.

Therefore, we design a bipartite graph representation for TPPs based on the structure of the objective function (\ref{eq_1}) and constraints (\ref{eq_2})-(\ref{eq_3}), with appropriate merging of the constraints and variables. Our bipartite graph representation is shown as Figure \ref{fig_4}. We represent the constraints associated only with \emph{products} (i.e., constraints (\ref{eq_2})) and the variables associated only with \emph{markets} / \emph{depot} (i.e.,  variables $v_i$ and $x_{ij}$) as two sets of nodes, and the constraints and variables associated with \emph{both products and markets} (i.e., constraints (\ref{eq_3}) and variables $z_{ik}$) as edges connecting the associated product nodes and market / depot nodes. Specifically, each variable $v_i$ is represented as a market / depot node, with the objective coefficients $c_{ij}$ of variables $x_{ij}$ represented by coordinate features of the market / depot nodes. Each constraint in (\ref{eq_2}) is represented as a product node, with the right-hand side values $d_k$ attached as node features. An edge connects a market node $i$ and a product node $k$ if the associated variable $z_{ik}$ is not always restricted to $0$ in constraints (\ref{eq_3}), that is, product $k$ is available at market $i$. The edge feature is a 2-dimensional vector, which consists of the corresponding objective coefficient $p_{ik}$ of variable $z_{ik}$ and the constraint coefficient $q_{ik}$ in constraints (\ref{eq_3}). For simplicity in the remainder of this paper, we use the term ``market nodes'' to refer to ``market / depot nodes'' if it does not cause confusion. \textcolor{black}{We remark that our bipartite graph representation can be viewed from a more intuitive perspective: it captures the two fundamental elements of TPPs---markets and products---as two sets of nodes, while their interdependencies, such as the available quantities and prices, are characterized using the edges between them.}

We emphasize that our bipartite graph representation effectively addresses the limitations of the complete or \emph{k-NN} graph representations discussed above. First, the dimensions of node and edge features are invariant to the number of markets and products, enabling the design of a size-agnostic policy network, such that it can be flexibly adapted to TPP instances of varying sizes. Second, the bipartite graph explicitly characterizes the relations between markets and products through edges connecting the market nodes and product nodes. The message-passing along edges can effectively capture the relational information between markets and products, which is important to the decision tasks.

\section{Policy Network} \label{Policy Network}

In this section, we introduce the architecture of our policy network, which is used for route construction in the first stage. In addition, at the end of this section, we will describe the basic training algorithm for updating the policy network during training.

\subsection{Policy Network Architecture}

\begin{figure*}[t]
\centering
\includegraphics[width=\textwidth, keepaspectratio]{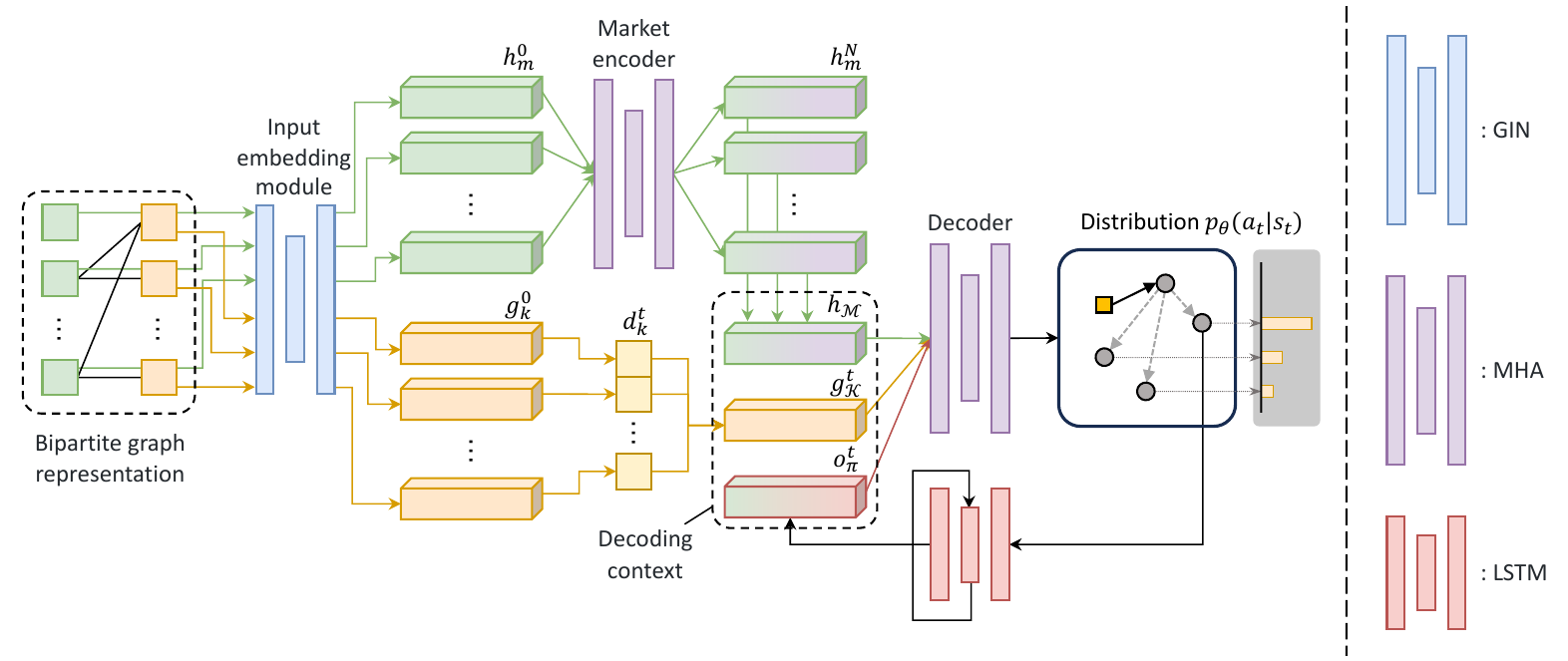}\\
\caption{The architecture of our policy network. The policy network is composed of 1) an input embedding module, 2) a market encoder, and 3) a decoder. The bipartite graph is embedded through the input embedding module and the market encoder. The decoder is executed iteratively: at step $t$, the decoder receives a decoding context, and outputs an action distribution $p_{\theta}(a_t \mid s_t)$, from which an action $a_t$ is selected, determining the next node to visit.}
\label{fig_2}
\end{figure*}

As described above, at each decision step $t$, the policy network, denoted as $\theta$, takes as input the state $s_t$ (consisting of the TPP instance $U$ being solved, which is represented as a bipartite graph, the partial route $\pi^t$, and the remaining demand $d^t$), and outputs a distribution $p_{\theta}(a_t | s_t)$ over the actions $a_t \in A_t$, determining the next node to visit. The architecture of our policy network is illustrated in Figure \ref{fig_2}. The policy network is composed of 1) an input embedding module, 2) a market encoder, and 3) a decoder.

The encoding process is designed to fully exploit the structure of bipartite graph. Specifically, the input embedding module takes the bipartite graph representation of instance $U$ as input, performing message-passing along edges to produce high-dimensional embeddings for the market and product nodes. These market node embeddings are then further processed through the market encoder to extract relevant information between each other. Following the encoding process, a decoder is executed iteratively to construct the route in a sequential manner. Specifically, at each decision step $t$, the decoder receives a decoding context that contains the embeddings of the bipartite graph and the dynamic part of current state $s_t$ (i.e., $\pi^t$ and $d^t$) as input. Based on this, the decoder outputs a distribution $p_{\theta}(a_t \mid s_t)$, i.e., the policy, from which an action $a_t$ is selected, determining the next node to visit. The dynamic part $\pi^t$ and $d^t$ of the state are updated accordingly, starting the next decoding step. Since the bipartite graph representation of $U$ is fixed throughout the state transitions in a complete MDP, we execute the input embedding module and market encoder only once at the initial state, and then the decoder is executed iteratively with each decision step of the MDP. In the following of this section, we provide a detailed introduction to the three components of our policy network.

\subsubsection{Input Embedding Module}

Taking the bipartite graph representation of instance $U$ as input, the input embedding module exploits the raw features and topology of the bipartite graph and produces a high-dimensional embedding for each market node and product node. GNN is an effective framework that naturally operates on graph-structured data for learning relational
information through message-passing between the nodes on the graph \citep{r_8}. Therefore, we incorporate GNNs into our input embedding module to embed the bipartite graph and capture the relations between markets and products.

\textcolor{black}{Considering the feature heterogeneity of the bipartite graph representation, we first normalize the node and edge features by the product demand values, and linearly project them into a uniform-dimensional space as initial embeddings.} To exploit the structure of the bipartite graph, we introduce a two-phase message-passing procedure that updates the initial node embeddings by aggregating information across the market nodes and product nodes: the first phase is performed to update the embeddings of product nodes, followed by the second phase that updates the embeddings of market nodes. Specifically, let $g_k^{\text{init}}$, $h_m^{\text{init}}$ and $e_{km}^{\text{init}}$ denote the initial embeddings of product node $k$, market node $m$, and edge $(k, m)$, respectively. The two-phase update proceeds as follows:
\begin{align*}
g_{k}^{0}=\operatorname{MLP}_{\mathcal{K}} & \left(\left(1+\epsilon \right) \cdot g_{k}^{\text{init}}+\sum_{m \in \mathcal{N}\left(k\right)} \operatorname{ReLU}\left(h_{m}^{\text{init}} + e_{mk}^{\text{init}} \right)\right), \\
h_{m}^{0}=\operatorname{MLP}_{\mathcal{M}} & \left(\left(1+\epsilon \right) \cdot h_{m}^{\text{init}}+\sum_{k \in \mathcal{N}\left(m\right)} \operatorname{ReLU}\left(g_{k}^{0} + e_{km}^{\text{init}} \right)\right),
\end{align*}
where $\operatorname{MLP}_{\mathcal{K}}$ and $\operatorname{MLP}_{\mathcal{M}}$ are multi-layer perceptrons (MLPs) for updating the product nodes and market nodes, respectively, and $\mathcal{N}\left(k\right)$ and $\mathcal{N}\left(m\right)$ denote the 1-hop neighbor node set of product node $k$ and market node $m$, respectively. Both phases can be seen as an extension of graph isomorphism network (GIN) \citep{r_9}, adapted to further take the bipartite structure into account and incorporate edge embeddings in the update. The updated product node embeddings $g_{k}^{0}$ for $k = 1, 2, \dots, |K|$ and market node embeddings $h_{m}^{0}$ for $m = 0, 1, \dots, |M|$ are used in the downstream market encoder and the decoder.

\textcolor{black}{We remark that our design behind the two-phase embedding procedure is intuitive. In the first phase, each product node gathers its supply information, including in which markets it is supplied, the coordinates of these markets, as well as the supply quantities and prices. This allows each product node to form an ``understanding'' of where and how it is supplied.
In the second phase, we make each market node collect higher-level information about its supplied products, including not only their total demands but also the information about how these products are supplied at other markets. This information is derived from the product node embeddings updated in the first phase. Consequently, the product node embedding $g_{k}^{0}$ aggregates its global supply information, and the market node embedding $h_{m}^{0}$ incorporates higher-level information such as the substitutive and complementary relations of product supply with other markets, which is important for making informed market selection decisions.}

\subsubsection{Market Encoder}

It is empirically observed from prior works for routing problems that deeper information extraction for the city nodes---i.e., market nodes in our task---can potentially improve the capability of the policy network \citep{r_36, r_53}.
Thus, we further process the market node embeddings through a market encoder for deeper information extraction.

Following the Transformer \citep{r_11} architecture, we process the market node embeddings through $N$ stacked attention layers. Each attention layer consists of two sublayers: an MHA layer that executes information aggregation for the market nodes and a node-wise MLP feed-forward layer to update the aggregated embeddings, where both sublayers add a residual connection \citep{r_15} and batch normalization (BN) \citep{r_12}:
\begin{align*}
& \hat{h}_{m}^{\ell+1}=\operatorname{BN}^{\ell}\left(h_{m}^{\ell}+\operatorname{MHA}_{m}^{\ell}\left(h_{0}^{\ell}, \ldots, h_{\left | M \right | }^{\ell}\right)\right), \\
& h_{m}^{\ell+1} = \operatorname{BN}^{\ell}\left(\hat{h}_{m}^{\ell+1} + \operatorname{MLP}^{\ell}\left(\hat{h}_{m}^{\ell+1}\right)\right),
\end{align*}
where $h_{m}^{\ell+1}$ is the embedding of market node $m$ after layer $\ell$ for $\ell=0, \ldots, N-1$. The MHA layer adopts a scaled dot-product attention mechanism, mapping a query vector and a set of key-value pairs to an output vector. For each attention head, a query vector $q_m^\ell$, a key vector $k_m^\ell$, and a value vector $v_m^\ell$ are calculated for each market node $m$ through a linear projection of its node embedding:
$$q_m^\ell=W_q^\ell h_{m}^\ell,\; k_m^\ell=W_k^\ell h_{m}^\ell,\; v_m^\ell=W_v^\ell h_{m}^\ell,$$
where $W_q^\ell$, $W_k^\ell$, and $W_v^\ell$ are trainable parameters. The attention output for market node $m$ is calculated as a weighted sum of the values $v_0^\ell, \ldots, v_{\left|M\right|}^\ell$, where the weight assigned to each value is the scaled dot-production of the query $q_m^\ell$ with the corresponding keys $k_0^\ell, \ldots, k_{\left|M\right|}^\ell$:
\begin{gather*}
u_{mn}^\ell = \frac{(q_m^\ell)^\mathrm{ T } \cdot k_n^\ell}{\sqrt{\text{dim}_k}},\\
\operatorname{Attention}^\ell(h_m^\ell) = \sum_{n=0}^{\left | M \right |}\operatorname{softmax}(u_{mn}^\ell) \cdot v_n^\ell,
\end{gather*}
where $\text{dim}_{k}$ denotes the dimension of key vectors. Multi-head attention can be interpreted as multiple attention functions in parallel, which enables the node to collect diverse messages from other nodes through different attention heads. The final MHA output for each market node $m$ is calculated by linearly projecting the concatenation of the attention heads:
$$\operatorname{MHA}_{m}^{\ell}\left(h_{0}^\ell, \ldots, h_{\left | M \right | }^\ell\right) = W_o^\ell \left[\text{head}_{m, 1}^\ell, \ldots, \text{head}_{m, H}^\ell\right],$$
where $\text{head}_{m, i}^\ell$ is an abbreviation for $\operatorname{Attention}^\ell_{i}(h_m^\ell)$ and $[\cdot, \cdot]$ denotes the concatenation operator.

After passing through $N$ attention layers, we calculate the mean of final market node embeddings as a global embedding $h_\mathcal{M}$, which is an aggregation of the final extracted information for the bipartite graph:
$$h_\mathcal{M} = \frac{1}{\left | M \right | + 1} \sum_{m=0}^{\left | M \right |} h_m^N.$$
\textcolor{black}{The global embedding $h_\mathcal{M}$ is then used as part of the decoding context, allowing the decoder to make well-informed decisions based on the global information of the instance.}

\subsubsection{Decoder}

The bipartite graph of the TPP instance $U$ is encoded through the input embedding module and market encoder. After that, the decoder is executed iteratively to construct the route, one node at a step. At each step $t=1, 2, ..., T$, the decoder takes as input a decoding context $h_d$ of the state $s_t$, and outputs a distribution $p_{\theta}(a_t \mid s_t)$ over the actions $a_t \in A_t$. Then, an action $a_t$ is selected based on the distribution. The corresponding market or depot is added to the end of the partial route, and the state is updated for the next decoding step, until a complete route is constructed.

The decoding context provides the decoder with the embedded information of current state $s_t$, including the instance $U$, the partial route $\pi^t$, and the remaining demand $d^t$. We form the decoding context $h_{d}$ as the concatenation of three parts: 1) the global embedding $h_\mathcal{M}$, 2) a demand context $g_\mathcal{K}^t$, and 3) a route context $o_\pi^t$:
$$h_{d} = [h_\mathcal{M}, g_\mathcal{K}^t, o_\pi^t].$$
First, as introduced above, the global embedding $h_\mathcal{M}$ is the aggregation of the ﬁnal extracted information from the bipartite graph, which contains global information of the instance $U$ being solved. Second, we form the demand context $g_\mathcal{K}^t$ as a weighted sum of product embeddings, where the weight for each product $k$ is its remaining demand $d_k^t$:
$$g_\mathcal{K}^t = \sum_{k=1}^{|K|} d_k^t \cdot g_{k}^{0}.$$
\textcolor{black}{The demand context $g_\mathcal{K}^t$ contains the information of remaining demand $d^t$, which helps the decoder consider which products are more inadequate when deciding the following markets.}

Third, the route context $o_\pi^t$ is the embedding for the partial route $\pi^t$, i.e., a sequence of market and depot nodes. We adopt an LSTM \citep{r_13} recurrent network, which updates the route context $o_\pi^t$ as follows:
$$o_\pi^t = \operatorname{LSTM}\left(o_\pi^{t-1}, h_{v^{t}}^N\right),$$
where $h_{v^{t}}^N$ is the final embedding of the node $v^t$ added to the partial route at the last step. \textcolor{black}{The LSTM helps capture the information about previously selected nodes and guiding the next decision in the route construction process.}

The decoder leverages the information from the decoding context $h_{d}$ to make decisions. It first processes the decoding context $h_{d}$ through a one-to-many attention layer, where the query is from the decoding context $h_d$, and the keys and values are from the final market / depot node embeddings $h_0^N, \ldots, h_{\left| M \right|}^N$. The action distribution $p_\theta \left(a_t \mid s_t \right)$ is computed as the single-head attention weights between the updated decoding context vector (denoted as $h_{d}^\prime$) and the final node embeddings $h_0^N, \ldots, h_{\left| M \right|}^N$ of markets / depot. We first computed attention scores between $h_{d}^\prime$ and the final market / depot node embeddings:
$$u_{dm} = \begin{cases}
 C \cdot \tanh\left( \frac{(q_d^\prime)^T \cdot k_m^N}{\sqrt{\text{dim}_k}} \right), \quad \text{if node } m \text{ is not masked},\\
-\infty, \, \quad \text{otherwise},
\end{cases}$$
where the query $q_d^\prime$ is from the updated decoding context vector $h_{d}^\prime$, and the keys $k_m^N$ are from the final market / depot node embeddings $h_m^N$ for $m=0,\ldots, |M|$. The attention scores are clipped using $\tanh$ function scaled by $C$, and the scores for the masked actions are set to $-\infty$. The final output distribution is computed using a $\operatorname{softmax}$ function:
$$p_\theta \left(a_t = m \mid s_t \right) = \frac{e^{u_{dm}}}{\sum_{n=0}^{\left | M \right |}e^{u_{dn}}}, \quad \text{for } m=0,\ldots, |M|.$$
Then, an action $a_t$ is selected based on the distribution, and the corresponding market / depot $m$ is added to the partial route.
Specifically, during training, the action is sampled from $p_\theta \left(a_t \mid s_t \right)$, while during inference, the action with the highest probability is selected greedily.

\begin{algorithm}[!b]
\SetAlgoLined
\caption{REINFORCE algorithm with greedy rollout baseline}
\label{alg_1}
\KwIn{Initial policy network parameter $\theta$, \\
\hspace{1.25cm} TPP instance distribution $\mathcal{P}$,\\
\hspace{1.25cm} number of epochs $E$, batch size $B$, \\
\hspace{1.25cm} steps per epoch $T$, learning rate $\epsilon$,\\
\hspace{1.25cm} significance $\alpha$.}
\KwOut{The learned policy network parameter $\theta$.}
Initialize baseline network parameter $\theta^{\text{BL}} \gets \theta$\;
\For{$e = 1, \ldots , E$}{
\For{$t= 1, \ldots , T$}{
Generate $B$ instances randomly from $\mathcal{P}$\;
Sample route $\pi_i \sim p_{\theta}(\pi_i | U_i),\; i \in\{1, \ldots, B\}$\;
Greedy rollout $b(U_i)$ from $p_{\theta^{\text{BL}}},\; i \in\{1, \ldots, B\}$ \tcp*{compute baseline}
$\nabla \mathcal{L} \leftarrow \sum_{i=1}^{B}\left(L\left({\pi}_{i} | U_i\right)-b(U_i)\right) \nabla_{\theta} \log p_{\theta}\left({\pi}_{i} | U_i\right)$ \tcp*{get gradient}
$\theta \gets \operatorname{Adam}\left(\theta, \nabla \mathcal{L} \right)$ \tcp*{update policy network}
}
\If{$\operatorname{OneSidedPairedTTest}\left(p_\theta, p_{\theta^{\text{BL}}}\right)<\alpha$}{
$\theta^{\text{BL}} \gets \theta$  \tcp*{update baseline network}
}
}
\end{algorithm}

\textcolor{black}{We highlight several advantages of our policy network design, which contribute to its capability, scalability, and efficiency. First, our policy network, based on GNNs and MHA, is agnostic to the size of instances being solved, given that the input feature dimensions are fixed by use of our bipartite graph representation. This enables the policy network to adapt directly to TPP instances with varying numbers of markets and products, without the need of parameter adjustment or retraining. Second, the embedding process for the static part $U$ and the dynamic part $\pi^t$ and $d^t$ of a state can be separated, which improves the computational efficiency of the policy network. Specifically, the embeddings of $U$ are pre-computed only once using the computationally intensive input embedding module and market encoder, while a computationally light decoder is executed repeatedly with varying decoding contexts that contain the dynamic information of states. This separation allows for efficient reuse of static embeddings and avoids redundant computations. Moreover, the forward-propagation of our policy network based on GNNs and MHA is highly parallelizable, which further improves the computational efficiency when used for solving TPPs.}

\subsection{Basic Training Algorithm} \label{Training Algorithm}

The policy network is updated based on the reward received after the solving process. As described earlier, given a TPP instance $U$, the policy network $\theta$ generates the final route $\pi$ with probability $p_{\theta}(\pi \mid U)=\prod_{t=1}^{T} p_{\theta}(a_t \mid s_t)$ during training, where the associated traveling cost can be directly accessed, and the purchasing cost is obtained by solving a transportation problem. To align with the notations in DRL, we let $L(\pi \mid U)$ denote the loss of $\pi$, which is defined as the negative reward for the terminal state. We update the policy network $\theta$ to minimize the expectation of loss $L(\pi \mid U)$, using the REINFORCE algorithm \citep{r_16} with baseline $b(U)$:
$$\nabla_{\theta} \mathcal{L}(\theta | U)=\mathbb{E}_{\pi \sim p_{\theta}(\pi | U)}\left[(L(\pi | U)-b(U)) \nabla_{\theta} \log p_{\theta}(\pi | U)\right],$$
where the baseline $b(U)$ can be seen as an estimation for the difficulty of instance $U$, which measures the relative advantage of route $\pi$ generated by the policy network. A well-defined baseline can significantly reduce the variance of gradients and accelerate training. We define the baseline $b(U)$ as the loss of the route obtained from a deterministic greedy rollout, where the action with the highest probability is always selected at each decision step, based on a baseline network $\theta^{\text{BL}}$ \citep{r_4}. Similar to the target Q-network in DQN \citep{r_39}, the baseline network $\theta^{\text{BL}}$ is a copy of policy network $\theta$, but fixed during each epoch to stabilize the baseline value. The baseline network $\theta^{\text{BL}}$ is periodically replaced by the latest policy network $\theta$ if the improvement of the policy is significant enough according to a paired t-test. Based on the gradient $\nabla_{\theta} \mathcal{L}(\theta \mid U)$, the parameters of policy network are updated using the Adam optimizer \citep{r_17}. The basic training algorithm introduced above is outlined as Algorithm \ref{alg_1}.

\section{Meta-learning Strategy} \label{Meta-learning Stategy}

Despite our proposed framework and the key DRL components, how to effectively learn a high-quality policy remains another challenge, especially for large-sized TPP instances.
First, due to the huge state-action space in solving TPPs, the training suffers from low sample efficiency. The policy network often makes random and suboptimal actions in the early training stage, requiring substantial trial-and-error before a good policy can be obtained. \textcolor{black}{This issue is exacerbated when training on large-sized TPP instances, where a randomly initialized policy network may fail to find a reasonable solution, resulting in potential training collapse (see Figure \ref{fig_3}(b)).} Second, while a well-trained policy can produce high-quality solutions for TPP instances from the same distribution as training, it tends to exhibit poor generalization performance on instances from different distributions. Generalization across different instance distributions is a highly desirable property for learning-based methods for solving CO problems.

In our approach, we address these limitations by pre-training an initialized policy network on a collection of varying instance distributions, such that it can take advantage of priorly learned cross-distribution knowledge and thus efficiently adapt to a new instance distribution using only a small amount of instances from that distribution for fine-tuning. Specifically, we propose a meta-learning strategy to learn such an initialized policy network. We consider a collection of instance distributions $\mathcal{D}_{\text{TPP}}=\left\{\mathcal{P}_i\right\}_{i=1}^D$, where each $\mathcal{P}_i$ defines a specific instance distribution (e.g., instances sharing the same number of markets and products). In our meta-training strategy, the training process comprises an outer-loop and several inner-loop optimizations at each iteration, corresponding to pre-training an initialized policy network and fine-tuning on a specific instance distribution, respectively. The outer-loop maintains a meta policy network $\theta$, which serves as the initialization $\theta_{\text{in}}^0$ for the inner-loop policy network. The inner-loop fine-tunes the policy network $\theta_{\text{in}}^0$ by performing $N$ updates on a specific instance distribution $\mathcal{P}_\text{in}$, which is drawn from $\mathcal{D}_{\text{TPP}}$ according to a sampling probability $p(\mathcal{D}_{\text{TPP}})=\left\{p(\mathcal{P}_i)\right\}_{i=1}^D$ (discrete uniform distribution in our implementation). The meta-objective is to learn a meta policy network $\theta$, which is a good initialization that can achieve low loss when adapting to a new instance distribution through $N$ fine-tuning updates:
$$\theta^{*}=\arg \min _{\theta} \mathbb{E}_{\mathcal{P}_\text{in} \sim p(\mathcal{D}_{\text{TPP}})} \mathbb{E}_{U \sim \mathcal{P}_\text{in}} \mathcal{L}\left(\theta_{\text{in}}^{N} \mid U \right),$$
where $\theta_{\text{in}}^{N}$ is the fine-tuned inner-loop policy network after $N$ updates from the initialization $\theta_{\text{in}}^0$.

\begin{algorithm}[!b]
\SetAlgoLined
\caption{Meta-training for the policy network}
\label{alg_2}
\KwIn{Initial meta policy network parameter $\theta$,\\
\hspace{1.25cm} TPP instance distributions set $\mathcal{D}_{\text{TPP}}$,\\
\hspace{1.25cm} number of epochs $E$, batch size $B$,\\
\hspace{1.25cm} outer steps per epoch $T$, outer step size $\beta$,\\
\hspace{1.25cm} inner steps per outer-loop $N$, learning rate $\epsilon$,\\
\hspace{1.25cm} significance $\alpha$.}
\KwOut{The trained meta policy network parameter $\theta$.}
Initialize baseline network parameter $\theta^{\text{BL}} \gets \theta$\;
\For{$e = 1, \ldots , E$}{
\For{$t = 1, \ldots , T$}{
Select a distribution $\mathcal{P}_{\text{in}} \sim p(\mathcal{D}_{\text{TPP}})$ \tcp*{start inner-loop}
Initialize inner model $\theta_{\text{in}} \gets \theta$\;
\For{$n = 1, \ldots , N$}{
Generate $B$ instances randomly from $\mathcal{P}_{\text{in}}$\;
Sample route $\pi_i \sim p_{\theta_{\text{in}}}(\pi_i | U_i),\; i \in\{1, \ldots, B\}$\;
Greedy rollout $b(U_i)$ from $p_{\theta^{\text{BL}}},\; i \in\{1, \ldots, B\}$ \tcp*{compute baseline}
$\nabla \mathcal{L} \leftarrow \sum_{i=1}^{B}\left(L\left({\pi}_{i} | U_i\right)-b(U_i)\right) \nabla_{\theta_{\text{in}}} \log p_{\theta_{\text{in}}}\left({\pi}_{i} | U_i\right)$ \tcp*{get gradient}
$\theta_{\text{in}} \gets \operatorname{Adam}\left(\theta_{\text{in}}, \nabla \mathcal{L} \right)$ \tcp*{inner-loop update}
}
$\theta \gets \theta + \beta \left(\theta_{\text{in}} - \theta\right)$  \tcp*{outer-loop update}
}
\If{$\operatorname{OneSidedPairedTTest}\left(p_\theta, p_{\theta^{\text{BL}}}\right)<\alpha$}{
$\theta^{\text{BL}} \gets \theta$\;
}
}
\end{algorithm}

In the inner-loop, we use the basic training algorithm (Section \ref{Training Algorithm}) to update the policy network $N$ times to obtain $\theta_{\text{in}}^{N}$ from $\theta_{\text{in}}^0$. To optimize the meta-objective, we validate the fine-tuned inner-loop policy network on a validation batch $\left\{U_i^{\prime}\right\}_{i=1}^B$, which is sampled from the same instance distribution $\mathcal{P}_\text{in}$ as the inner-loop training. The meta-gradient to update the meta policy network $\theta$ is according to the gradient chain rule as follows:
$$\nabla_{\theta} \mathcal{L} \left(\theta_{\text{in}}^{N}\right) = \frac{1}{B} \sum_{i=1}^{B} \nabla_{\theta_{\text{in}}^{N}} \mathcal{L}\left(\theta_{\text{in}}^{N} \mid U_i^{\prime}\right) \cdot \frac{\partial \theta_{\text{in}}^{N}}{\partial \theta_{\text{in}}^0}.$$
Unfortunately, the meta-gradient involves a second-order derivative in $\frac{\partial \theta_{\text{in}}^{N}}{\partial \theta_{\text{in}}^0}$, which is significantly expensive in computation. Therefore, we instead employ a computationally efficient first-order approximation for the meta-gradient. The update for the meta policy network is approximated as moving towards the parameters of the final inner-loop policy network \citep{r_19}:
$$\theta \gets \theta + \beta \left(\theta_{\text{in}}^{N} - \theta\right),$$
where $\beta$ is the outer-loop step size. The procedure of our meta-learning strategy is outlined as Algorithm \ref{alg_2}.

Through meta-training, the meta policy network can serve as an effective initialization that achieves good performance on a new instance distribution through a few fine-tuning steps. This is especially beneficial for large-sized instances, which tend to experience training collapse if the policy network is trained from scratch. Moreover, the meta policy network can learn cross-distribution knowledge through meta-learning, enhancing its ability to generalize across varied instance sizes and distributions, even showing good zero-shot generalization performance on instances from a different distribution that is never seen during training.

\section{Experiments} \label{Experiments}

We conduct extensive experiments on various synthetic TPP instances and the TPPLIB benchmark to evaluate our DRL-based approach, comparing it against several well-established TPP heuristics. Furthermore, we validate the zero-shot generalization ability of our policy network on larger-sized TPP instances that are unseen during training. In addition, an ablation study to analyze the contributions of each component in our framework and extended results are presented in the appendices.

\subsection{Experimental Settings}

\subsubsection{Instances}

We first randomly generate synthetic TPP instances of sizes $\left(|M|, |K|\right) = (50, 50)$, $(50, 100)$, $(100, 50)$, and $(100, 100)$ for training and evaluation, where $|M|$ is the number of markets and $|K|$ is the number of products. The generation of U-TPP and R-TPP instances follows the standard procedure introduced in \citep{r_22}, corresponding to Class 3 and Class 4 in the TPPLIB benchmark \citep{r_35}. For the U-TPP instance, $(|M|+1)$ integer coordinates (including the depot) are randomly generated within the $[0, 1000] \times [0, 1000]$ square according to a uniform distribution, and the traveling costs are defined as the truncated Euclidean distances. Each product $k \in K$ is supplied at $|M_k|$ randomly selected markets, where $|M_k|$ is uniformly generated in $\left[1, |M|\right]$. The price $p_{ik}$ of each product $k \in K$ at each market $i \in M_k$ is randomly sampled in $[1, 10]$. The R-TPP instances are generated in a similar manner as U-TPP, with an additional limit $q_{ik}$ on the available supply quantities, which is randomly sampled in $[1, 15]$ for each product $k \in K$ and each market $i \in M_k$. A parameter $\lambda$ is used to control the number of markets in a feasible solution through the product demand $d_{k}:=$ $\left\lceil\lambda \max _{i \in M_{k}} q_{ik}+(1-\lambda) \sum_{i \in M_{k}} q_{ik}\right\rceil, 0<\lambda<1$, for $k \in K$. We choose $\lambda=\left\{0.99, 0.95, 0.9\right\}$ in our experiments, corresponding to the most difficult instances in TPPLIB \citep{r_1}.

We define an instance distribution as U-TPP instances sharing the same $\left(|M|, |K|\right)$, or R-TPP instances sharing the same $(|M|, |K|, \lambda)$. For each instance distribution, training instances are generated on-the-fly, and a collection of 1000 instances is generated for evaluation. Furthermore, we generate larger-sized instances for evaluating the zero-shot generalization ability of the policy network trained on smaller instances. In addition, the well-known TPPLIB benchmark is also used for evaluation, which includes 5 instances for each instance distribution, called a category.

\subsubsection{Baseline Methods}

We select several well-established TPP heuristics reported in \citep{r_1}, including GSH, CAH, and TRH. A widely used and effective practice is to incorporate constructive heuristics (e.g., GSH and CAH) with TRH to remove redundant markets as soon as a solution is produced. Additionally, the solution can potentially be further improved by using a TSP solver to re-sequence the visited markets. In our experiments, we implement GSH and CAH, both followed by TRH and the TSP re-sequence, as the baselines, denoted as ``GSH+TRH'' and ``CAH+TRH'', respectively. 

Note that we do not include exact methods in the evaluation on synthetic instances, as it can take several days to weeks to solve an evaluation instance set using exact methods (even for U-TPP). For a comparison of the optimality gap, we reference the best-known solutions in the literature \citep{r_35} on the TPPLIB benchmark instances to compare the optimality gaps of the baseline heuristics and our DRL-based approach.

\begin{table}[t]
\fontsize{11}{11}\selectfont
\renewcommand\arraystretch{1.1}
\caption{Training configuration} \label{Training configuration}
\vspace{5pt}
\centering
\begin{tabular}{lcr}
\toprule
Hyper-parameter & \hspace{60pt} & Value \\
\midrule
Embedding dimension & \hspace{10pt} & 128 \\
Num of market encoder layers & \hspace{10pt} & 3 \\
Key vector dimension & \hspace{10pt} & 16 \\
Num of attention heads & \hspace{10pt} & 8 \\
Clipping scale for $\tanh$ & \hspace{10pt} & 10.0 \\
\midrule
Num of epochs & \hspace{10pt} & 100 \\
Batch size & \hspace{10pt} & 512 \\
Steps per epoch & \hspace{10pt} & 2500 \\
Learning rate & \hspace{10pt} & 1e-4 \\
T-test significance & \hspace{10pt} & 0.05 \\
\midrule
Outer steps per epoch & \hspace{10pt} & 2500 \\
Inner steps per outer-loop & \hspace{10pt} & 2 \\
Outer step size & \hspace{10pt} & 0.8 \\
\bottomrule
\end{tabular}
\end{table}

\subsubsection{Training and Inference of the Policy Network}

The detailed training configuration for the policy network is presented in Table \ref{Training configuration}. Both training and inference are performed on a machine equipped with an AMD EPYC 7742 CPU at 2.3 GHz and a single GeForce RTX-4090 GPU.

For evaluating our DRL-based approach, we perform inference of the learned policy network using greedy decoding (similar to the greedy rollout baseline) with instance augmentation \citep{r_36}. The produced route and purchasing plan can directly serve as an end-to-end solution to the TPP instance, which we denote as ``RL-E2E''. Similar to the baseline heuristics, we can also apply a post-optimization procedure (TRH and TSP re-sequence) to further improve the end-to-end solution, denoted as ``RL+TRH''. For a fair comparison, the baseline heuristics and our approach are implemented on the same machine and environment.

\subsection{Training Performance} \label{Training Performance}

First, we present the training performance of the policy network on different instance distributions, highlighting the necessity and effectiveness of our meta-learning strategy.

\begin{figure}[t]
\centering
\includegraphics[width=3.7in, keepaspectratio]{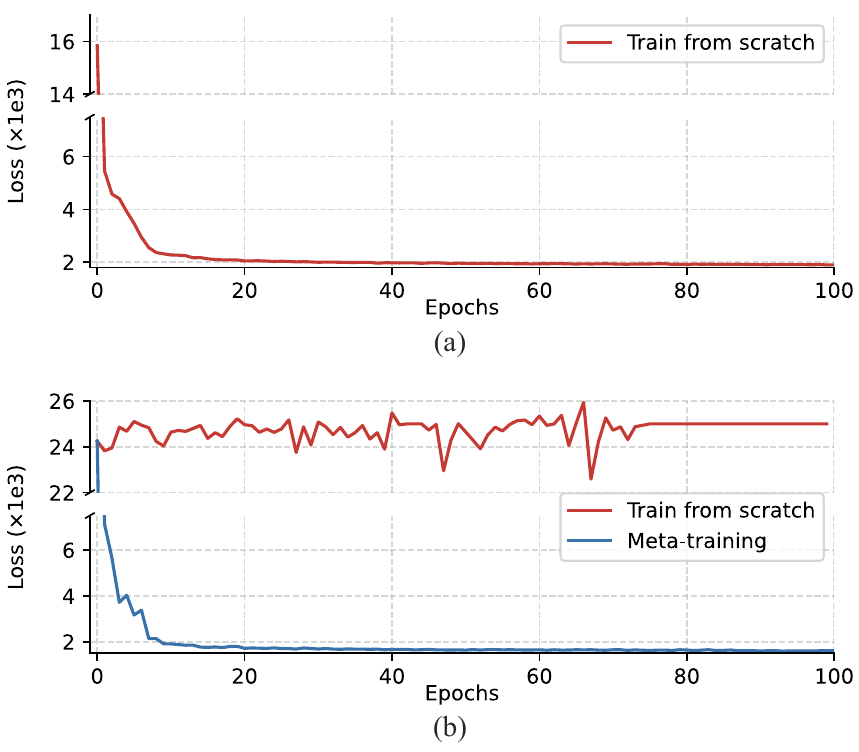}\\
\caption{The training curve on (a) U-TPP instances of size $\left(|M|, |K|\right) = (50, 50)$, and (b) U-TPP instances of size $\left(|M|, |K|\right) = (100, 50)$.}
\label{fig_3}
\end{figure}

The policy network is trained from scratch for U-TPP and R-TPP instances of sizes $\left(|M|, |K|\right) = (50, 50)$ and $(50, 100)$. For brevity, we illustrate the training curve on U-TPP instances of size $\left(|M|, |K|\right) = (50, 50)$ in Figure \ref{fig_3}(a), which is plotted based on the average loss on the evaluation instance set. The stable convergence demonstrates that our policy network is capable of learning an effective policy for small-sized TPPs when trained from scratch using the basic training algorithm.

However, the policy network suffers serious training collapse when trained from scratch on large-sized U-TPP and R-TPP instances of sizes $\left(|M|, |K|\right) = (100, 50)$ and $(100, 100)$. As shown in Figure \ref{fig_3}(b), the policy network trained from scratch on U-TPP instances of size $\left(|M|, |K|\right) = (100, 50)$ is stuck in a high-loss region, and the loss completely fails to decrease throughout training. This is mainly because the state-action space is so large that the policy network cannot find a reasonable solution if the training starts from a random initialization. Therefore, we instead apply our meta-learning strategy to train the policy network for U-TPP and R-TPP instances with $\left(|M|, |K|\right) = (100, 50)$ and $(100, 100)$. As shown in Figure \ref{fig_3}(b), the utilization of meta-learning strategy significantly promotes the training for an effective policy.

\subsection{Results on Synthetic Instances}

Following the training procedures described above, we train the policy network from scratch for U-TPP and R-TPP instances of sizes $\left(|M|, |K|\right) = (50, 50)$ and $(50, 100)$, and apply the meta-learning strategy for training U-TPP and R-TPP instances of sizes $\left(|M|, |K|\right) = (100, 50)$ and $(100, 100)$. Then, the learned policy network is evaluated on the synthetic evaluation instance set for each instance distribution. \textcolor{black}{We report the metrics of 1) the average objective value of the obtained solution and 2) the average runtime to find this solution (in seconds) per instance.} The results are reported in Table \ref{synthetic U-TPP} for U-TPP and Table \ref{synthetic R-TPP} for R-TPP.

In terms of the objective value, the performance of baseline heuristics varies across different instance distributions. In contrast, our DRL-based approach consistently outperforms both GSH+TRH and CAH+TRH on all synthetic U-TPP and R-TPP instance sets of different sizes and distributions. Notably, the end-to-end solutions (RL-E2E) can already achieve better objective value than the hand-crafted baseline heuristics, and the post-optimization procedure (RL+TRH) further improves the solution, with only a slight increase in runtime. The results demonstrate that our DRL-based approach can learn effective policies for solving TPPs of different sizes and distributions. 

\begin{table}[!b]
\setlength{\tabcolsep}{8pt}
\fontsize{8}{8}\selectfont
\renewcommand\arraystretch{1.2}
\caption{Results on synthetic U-TPP instances} \label{synthetic U-TPP}
\vspace{5pt}
\centering
\begin{tabular}{cc|cccccccc}
\toprule
\multicolumn{2}{c}{Instance} & \multicolumn{2}{c}{GSH + TRH} & \multicolumn{2}{c}{CAH + TRH} & \multicolumn{2}{c}{RL - E2E} & \multicolumn{2}{c}{RL + TRH} \\
\cmidrule(lr){1-2} \cmidrule(lr){3-4} \cmidrule(lr){5-6} \cmidrule(lr){7-8} \cmidrule(lr){9-10}
$|M|$ & $|K|$ & Obj. & Time & Obj. & Time & Obj. & Time & Obj. & Time \\
\midrule
50 & 50   & 2221 & 0.006 & 1910 & 0.017 & 1897 & 0.017 & \textbf{1857} & 0.024 \\
50 & 100  & 2750 & 0.008 & 2552 & 0.033 & 2542 & 0.025 & \textbf{2446} & 0.033 \\
100 & 50  & 2050 & 0.011 & 1571 & 0.033 & 1563 & 0.020 & \textbf{1524} & 0.027 \\
100 & 100 & 2542 & 0.016 & 2185 & 0.072 & 2111 & 0.027 & \textbf{2044} & 0.036 \\
\bottomrule
\end{tabular}
\end{table}

\begin{table}[!b]
\setlength{\tabcolsep}{8pt}
\fontsize{8}{8}\selectfont
\renewcommand\arraystretch{1.2}
\caption{Results on synthetic R-TPP instances}  \label{synthetic R-TPP}
\vspace{5pt}
\centering
\begin{tabular}{ccc|cccccccc}
\toprule
\multicolumn{3}{c}{Instance} & \multicolumn{2}{c}{GSH + TRH} & \multicolumn{2}{c}{CAH + TRH} & \multicolumn{2}{c}{RL - E2E} & \multicolumn{2}{c}{RL + TRH} \\
\cmidrule(lr){1-3} \cmidrule(lr){4-5} \cmidrule(lr){6-7} \cmidrule(lr){8-9} \cmidrule(lr){10-11}
$|M|$ & $|K|$ & $\lambda$ & Obj. & Time & Obj. & Time & Obj. & Time & Obj. & Time \\
\midrule
50 & 50 & 0.99 & 2152 & 0.016 & 2257 & 0.073 & 2032 & 0.020 & \textbf{1954} & 0.030 \\
50 & 100 & 0.99 & 2671 & 0.032 & 2863 & 0.161 & 2567 & 0.026 & \textbf{2466} & 0.042 \\
100 & 50 & 0.99 & 2062 & 0.058 & 2174 & 0.243 & 1753 & 0.029 & \textbf{1711} & 0.043 \\
100 & 100 & 0.99 & 2578 & 0.142 & 2853 & 0.704 & 2302 & 0.033 & \textbf{2235} & 0.053 \\
\midrule
50 & 50 & 0.95 & 2845 & 0.044 & 2914 & 0.124 & 2645 & 0.025 & \textbf{2594} & 0.036 \\
50 & 100 & 0.95 & 3569 & 0.095 & 3709 & 0.243 & 3457 & 0.034 & \textbf{3368} & 0.059 \\
100 & 50 & 0.95 & 3384 & 0.300 & 3440 & 0.598 & 3027 & 0.040 & \textbf{2980} & 0.073 \\
100 & 100 & 0.95 & 4281 & 0.696 & 4405 & 1.499 & 3993 & 0.053 & \textbf{3920} & 0.111 \\
\midrule
50 & 50 & 0.9 & 3910 & 0.099 & 3965 & 0.201 & 3704 & 0.033 & \textbf{3644} & 0.052 \\
50 & 100 & 0.9 & 5080 & 0.195 & 5137 & 0.350 & 4927 & 0.043 & \textbf{4855} & 0.080 \\
100 & 50 & 0.9 & 5293 & 0.789 & 5310 & 1.218 & 4956 & 0.060 & \textbf{4900} & 0.124 \\
100 & 100 & 0.9 & 6963 & 1.666 & 7014 & 2.637 & 6672 & 0.073 & \textbf{6660} & 0.186 \\
\bottomrule
\end{tabular}
\end{table}

\begin{figure*}[!b]
\centering
\includegraphics[width=6.1in, keepaspectratio]{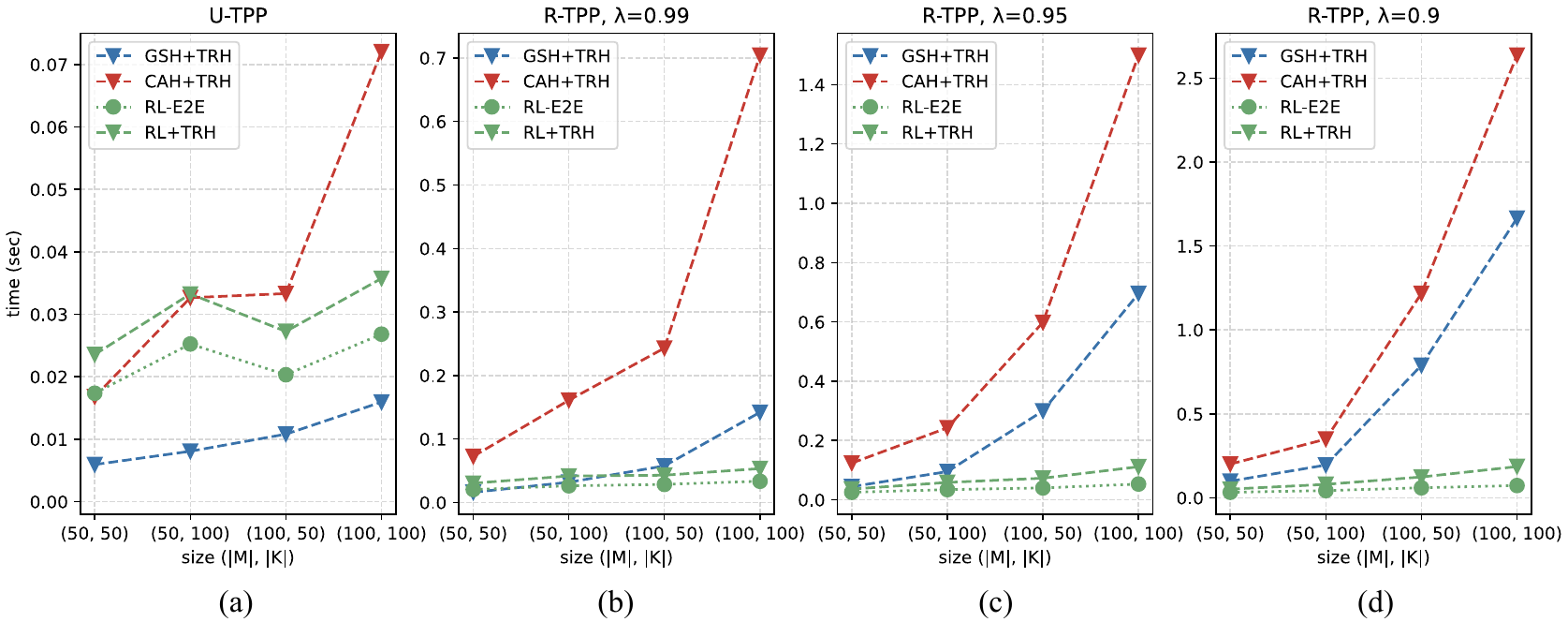}\\
\caption{The runtime of the baseline methods and our DRL-based approach on the synthetic (a) U-TPP, (b) R-TPP with $\lambda=0.99$, (c) R-TPP with $\lambda=0.95$, and (d) R-TPP with $\lambda=0.9$ instances.}
\label{fig_5}
\end{figure*}

In terms of the runtime, all methods, including the baseline heuristics and our DRL-based approach, can produce a solution within an average of 0.1 seconds for U-TPP instances of all four sizes. Generally, GSH+TRH is slightly faster, but the difference is approximately negligible in practice. However, for more challenging R-TPP instances, the runtime differences become more pronounced. As shown in Figure \ref{fig_5}(b)-(d), the runtime of baseline heuristics increases significantly with the instance size. GSH+TRH and CAH+TRH take an average of 1.666 seconds and 2.637 seconds, respectively, to produce a solution for R-TPP instances with $|M|=100, |K|=100, \lambda=0.9$. This is mainly because they need to frequently compute the objective values and savings each time a new market is inserted into the route. In contrast, the actions in our DRL-based approach are entirely based on the forward-propagation of the policy network. Therefore, the solution time is approximately linear with the number of visited markets in the constructed route. On synthetic R-TPP instances, the average solution time of RL-E2E is consistently below 0.1 seconds, and the average solution time of RL+TRH is consistently below 0.2 seconds.

\subsection{Results on TPPLIB Benchmark Instances}

We next evaluate our approach on the TPPLIB benchmark, by directly running the policy network (trained on synthetic instances) on the benchmark instances. The TPPLIB benchmark contains 5 instances for each instance distribution as a category. The results are summarized in Table \ref{TPPLIB benchmark}. The first block of the table (labeled with ``EEuclideo'') reports the results on U-TPP instances, and the following three blocks (labeled with ``CapEuclideo'') report the results on R-TPP instances with $\lambda=0.99, 0.95, 0.9$, respectively. Each row presents the average result for a category of 5 benchmark instances. For example, ``CapEuclideo.50.100.95'' refers to R-TPP instances with $|M|=50, |K|=100, \lambda=0.95$. We report 1) the average objective value, 2) the average optimality gap, and 3) the average runtime (in seconds) for each category of the TPPLIB benchmark.

\begin{table}[!b]
\setlength{\tabcolsep}{3pt}
\fontsize{8}{8}\selectfont
\renewcommand\arraystretch{1.2}
\caption{Results on TPPLIB benchmark instances} \label{TPPLIB benchmark}
\vspace{5pt}
\centering
\begin{tabular}{c|cccccccccccccc}
\toprule
\multirow{2.5}{*}{Instance} & \multicolumn{2}{c}{Opt.} & \multicolumn{3}{c}{GSH + TRH} & \multicolumn{3}{c}{CAH + TRH} & \multicolumn{3}{c}{RL - E2E} & \multicolumn{3}{c}{RL + TRH} \\
\cmidrule(lr){2-3} \cmidrule(lr){4-6} \cmidrule(lr){7-9} \cmidrule(lr){10-12} \cmidrule(lr){13-15}
 & Obj. & Time & Obj. & Gap & Time & Obj. & Gap & Time & Obj. & Gap & Time & Obj. & Gap & Time  \\
\midrule
EEuclideo.50.50 \, \, & 1482 & 4 & 1779 & 17.20\% & 0.007 & 1643 & 9.97\% & 0.015 & 1524 & 2.63\% & 0.014 & \textbf{1497} & \textbf{1.06\%} & 0.019 \\
EEuclideo.50.100 \, & 2417 & 6 & 2652 & 9.74\% & 0.008 & 2726 & 13.82\% & 0.032 & 2588 & 7.01\% & 0.025 & \textbf{2446} & \textbf{1.03\%} & 0.037 \\
EEuclideo.100.50 \, & 1655 & 72 & 1979 & 24.89\% & 0.009 & 1796 & 9.45\% & 0.031 & 1701 & 2.96\% & 0.019 & \textbf{1688} & \textbf{2.30\%} & 0.027 \\
EEuclideo.100.100 & 2085 & 183 & 2388 & 15.47\% & 0.013 & 2251 & 8.65\% & 0.075 & 2220 & 7.09\% & 0.025 & \textbf{2146} & \textbf{2.99\%} & 0.064 \\
\midrule
CapEuclideo.50.50.99 \, \, & 1862 & 6 & 2189 & 20.86\% & 0.017 & 2243 & 23.72\% & 0.089 & 2047 & 9.09\% & 0.020 & \textbf{1929} & \textbf{3.51\%} & 0.027 \\
CapEuclideo.50.100.99 \, & 2313 & 7 & 2578 & 12.30\% & 0.031 & 2710 & 18.23\% & 0.160 & 2483 & 7.18\% & 0.025 & \textbf{2394} & \textbf{3.33\%} & 0.033 \\
CapEuclideo.100.50.99 \, & 1504 & 58 & 1951 & 29.97\% & 0.061 & 1988 & 35.73\% & 0.179 & 1561 & 3.91\% & 0.025 & \textbf{1531} & \textbf{1.84\%} & 0.034 \\
CapEuclideo.100.100.99 & 1865 & 134 & 2283 & 21.76\% & 0.118 & 2406 & 27.95\% & 0.686 & 1955 & 4.86\% & 0.028 & \textbf{1914} & \textbf{2.79\%} & 0.039 \\
\midrule
CapEuclideo.50.50.95 \, \, & 2444 & 10 & 2904 & 21.16\% & 0.036 & 2751 & 15.61\% & 0.104 & 2643 & 7.40\% & 0.028 & \textbf{2581} & \textbf{5.09\%} & 0.040 \\
CapEuclideo.50.100.95 \, & 3187 & 23 & 3441 & 7.97\% & 0.072 & 3672 & 15.75\% & 0.234 & 3421 & 7.35\% & 0.036 & \textbf{3299} & \textbf{3.56\%} & 0.054 \\
CapEuclideo.100.50.95 \, & 2860 & 466 & 3144 & 9.85\% & 0.199 & 3221 & 12.73\% & 0.629 & 3026 & 5.93\% & 0.039 & \textbf{2962} & \textbf{3.66\%} & 0.063 \\
CapEuclideo.100.100.95 & 3555 & 1178 & 3991 & 12.31\% & 0.521 & 4096 & 15.24\% & 1.249 & 3769 & 5.94\% & 0.049 & \textbf{3664} & \textbf{3.00\%} & 0.094 \\
\midrule
CapEuclideo.50.50.9 \, \, & 3571 & 28 & 3927 & 10.05\% & 0.061 & 3873 & 8.49\% & 0.145 & 3744 & 4.73\% & 0.036 & \textbf{3673} & \textbf{2.81\%} & 0.053 \\
CapEuclideo.50.100.9 \, & 4668 & 30 & 4961 & 6.33\% & 0.128 & 5046 & 8.07\% & 0.289 & 4876 & 4.47\% & 0.043 & \textbf{4834} & \textbf{3.57\%} & 0.067 \\
CapEuclideo.100.50.9 \, & 4674 & 243 & 4981 & 6.64\% & 0.439 & 5106 & 9.26\% & 0.999 & 4891 & 4.66\% & 0.053 & \textbf{4825} & \textbf{3.25\%} & 0.097 \\
CapEuclideo.100.100.9 & 6442 & 537 & 6961 & 8.11\% & 1.668 & 6850 & 6.43\% & 2.307 & 6637 & 3.01\% & 0.070 & \textbf{6534} & \textbf{1.42\%} & 0.156 \\
\midrule
\textcolor{black}{Average} & \textcolor{black}{2912} & \textcolor{black}{187} & \textcolor{black}{3257} & \textcolor{black}{14.66\%} & \textcolor{black}{0.212} & \textcolor{black}{3274} & \textcolor{black}{14.94\%} & \textcolor{black}{0.451} & \textcolor{black}{3068} & \textcolor{black}{5.51\%} & \textcolor{black}{0.033} & \textcolor{black}{\textbf{2995}} & \textcolor{black}{\textbf{2.83\%}} & \textcolor{black}{0.057} \\
\bottomrule
\end{tabular}
\end{table}

\begin{figure*}[!t]
\centering
\includegraphics[width=6.1in, keepaspectratio]{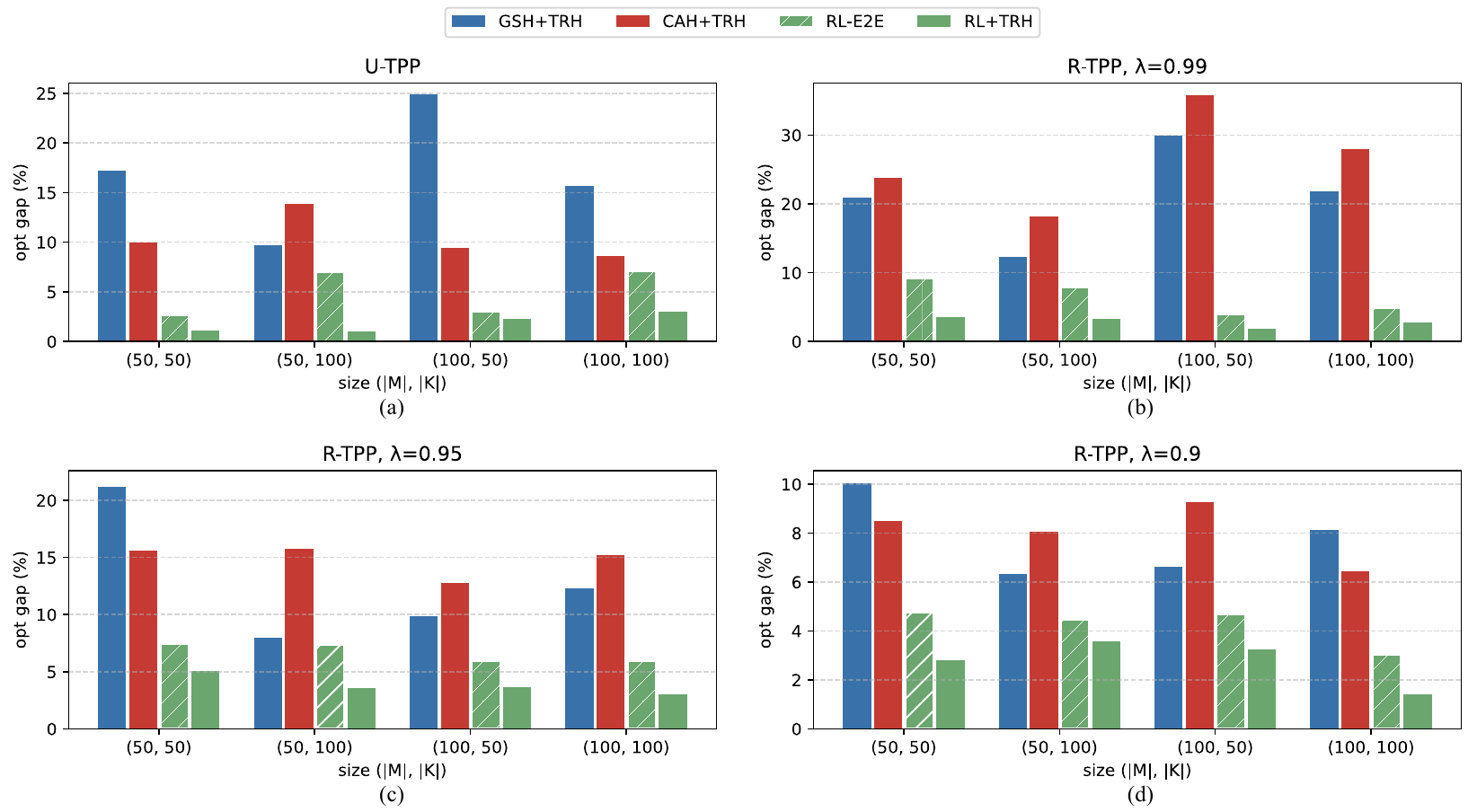}\\
\caption{The optimality gaps of the baseline methods and our DRL-based approach on (a) U-TPP, (b) R-TPP with $\lambda=0.99$, (c) R-TPP with $\lambda=0.95$, and (d) R-TPP with $\lambda=0.9$ in the TPPLIB benchmark.}
\label{fig_6}
\end{figure*}

The results demonstrate that the policy network, trained on synthetic instances, can be effectively applied to TPPLIB benchmark instances from the same instance distribution. Our DRL-based approach outperforms the baseline heuristics by a large margin in terms of the optimality gap, especially for large-sized R-TPP instances. In detail, the average optimality gaps of solutions produced using our RL+TRH method are consistently within 6\% on each category of the TPPLIB benchmark in the experiment, yielding a reduction ranging from 40\% to 90\% compared to the baseline heuristics. The comparison of optimality gaps is illustrated in Figure \ref{fig_6}. Besides, similar to the results on synthetic instances, our DRL-based approach also shows a significant advantage in runtime.

\subsection{Generalization Performance on Large-Sized Instances}

Furthermore, we validate the zero-shot generalization performance of our DRL-based approach. We employ the meta-learning strategy to train the policy network on TPP instances of sizes $\left(|M|, |K|\right) = (50, 50)$, $(50, 100)$, $(100, 50)$, and $(100, 100)$, and then directly apply the meta policy network to solve larger-sized problem instances that are never seen in the training stage, without any fine-tuning. The results are summarized in Table \ref{large instance}. The first block of Table \ref{large instance} presents the results on U-TPP instances of sizes $\left(|M|, |K|\right)=(150, 150)$, $(200, 200)$, and $(300, 300)$, and the second block presents the results on R-TPP instances with $|M|=150, |K|=150$ and $\lambda=0.99, 0.95, 0.9$.

It is demonstrated that our DRL-based approach still shows an advantage over the baseline heuristics when generalizing to U-TPP instances of sizes $\left(|M|, |K|\right) = (150, 150)$ and $(200, 200)$, while the largest instances for training are only of size $\left(|M|, |K|\right) = (100, 100)$. However, we observe a drop in the zero-shot generalization performance on U-TPP instances of size over $\left(|M|, |K|\right) = (300, 300)$ and R-TPP instances of size over $\left(|M|, |K|\right) = (150, 150)$. In future work, we shall further explore the generalization and more efficient learning techniques for larger-sized TPP instances.

\begin{table}[!b]
\setlength{\tabcolsep}{8pt}
\fontsize{8}{8}\selectfont
\renewcommand\arraystretch{1.2}
\caption{Results on larger-sized instances} \label{large instance}
\vspace{5pt}
\centering
\begin{tabular}{ccc|cccccccc}
\toprule
\multicolumn{3}{c}{Instance} & \multicolumn{2}{c}{GSH + TRH} & \multicolumn{2}{c}{CAH + TRH} & \multicolumn{2}{c}{RL - E2E} & \multicolumn{2}{c}{RL + TRH} \\
\cmidrule(lr){1-3} \cmidrule(lr){4-5} \cmidrule(lr){6-7} \cmidrule(lr){8-9} \cmidrule(lr){10-11}
$|M|$ & $|K|$ & $\lambda$ & Obj. & Time & Obj. & Time & Obj. & Time & Obj. & Time \\
\midrule
150 & 150 & / & 2672 & 0.027 & 2460 & 0.190 & 2514 & 0.028 & \textbf{2380} & 0.075 \\
200 & 200 & / & 2885 & 0.034 & 2454 & 0.259 & 2485 & 0.040 & \textbf{2262} & 0.142 \\
300 & 300 & / & 3445 & 0.041 & \textbf{3074} & 0.373 & 3206 & 0.045 & 3088 & 0.198 \\
\midrule
150 & 150 & 0.99 & 2918 & 0.265 & 3056 & 1.202 & 2576 & 0.039 & \textbf{2519} & 0.063 \\
150 & 150 & 0.95 & \textbf{5642} & 1.318 & 5884 & 1.862 & 6133 & 0.074 & 6014 & 0.150 \\
150 & 150 & 0.9 & \textbf{10008} & 3.027 & 10199 & 4.908 & 10707 & 0.095 & 10155 & 0.205 \\
\bottomrule
\end{tabular}
\end{table}

\vspace{15pt}

In addition, to get a better understanding of the contributions of components in our framework, we conduct a series of ablation studies. The results are presented in \ref{Ablation Study}. We also present extended results that include comparisons with state-of-the-art neural VRP solvers and a validation of the cross-problem compatibility with other VRP variants in \ref{Extended Results} to further highlight the effectiveness and compatibility of our approach.

\section{Conclusion} \label{Conclusion}

In this paper, we have presented a novel DRL-based approach for solving TPPs, which exploits the idea of ``\emph{solve separately, learn globally}''. Namely, we break the solution task into two separate stages at the operational level: route construction and purchase planning, while learning a policy network towards optimizing the global solution objective. Built on this framework, we proposed a bipartite graph representation for TPPs and designed a policy network that effectively extracts information from the bipartite graph for route construction. Moreover, we introduced a meta-learning strategy, which significantly enhances the training stability and efficiency on large-sized instances and improves the generalization ability. Experimental results demonstrate that our DRL-based approach can significantly outperform well-established TPP heuristics in both solution quality and runtime, while also showing good generalization performance. Future works shall further explore better and more efficient generalization to larger-sized TPP instances and a more generic framework to address other TPP variants, such as multi-vehicle TPP, dynamic TPP, etc.


\addcontentsline{toc}{section}{\numberline{}References}
\bibliographystyle{elsarticle-harv} 
\bibliography{reference}

\newpage

\appendix
\renewcommand{\thetable}{\arabic{table}}

\section{Related Work} \label{Related Work}
Here we review the existing methods for TPPs, DRL-based methods for routing problems, and meta-learning.

\subsection*{A.1. Traveling Purchaser Problem}

The TPP is an important combinatorial optimization problem in logistics and manufacturing, with broad applications such as purchasing required raw materials for manufacturing factories \citep{r_23}, or scheduling a set of jobs for certain machines \citep{r_49}, etc. Notably, the TSP is a special case of TPPs \citep{r_1}.

Existing methods for solving TPPs can be categorized into exact methods and heuristics \citep{r_1}. The first exact method, introduced in \citep{r_57}, solved U-TPP using dynamic programming. Subsequently, \cite{r_58} developed a branch-and-bound algorithm, capable of solving TPP instances with up to 20 markets and 100 products, or 25 markets and 50 products. \cite{r_22} further enhanced this approach through a branch-and-cut algorithm that exploits dynamic generation of variables and separation of constraints. \cite{r_41} extended this methodology to general TPPs. Despite a significant advance in the size of TPPs solved to optimality, their computational cost remains unaffordable in real-time applications.

Alternatively, various heuristic methods have emerged to generate high-quality solutions within a reasonable time. GSH \citep{r_24}, the first constructive heuristic for TPPs, starts from the depot and iteratively inserts the unvisited market offering the largest cost saving. \cite{r_25} introduced TRH, a reduction-based heuristic that starts with an initial route containing a subset of markets, and iteratively drops the market yielding the largest cost reduction. The CAH \citep{r_26} adopts a product-focused approach, constructing a least-cost solution for the first product and, in subsequent iterations, adding each product to the solution in a least-cost manner. In practice, GSH and CAH are usually followed by TRH and a TSP solver to remove redundant markets and optimize market sequencing for further improvements. However, these heuristics rely heavily on substantial specialized expertise and knowledge, and their performance remains limited. Additional works, such as local search methods, explore improved solutions at the expense of extra computational time \citep{r_31, r_30}, but, overall, there was no significant breakthrough in heuristics for TPPs.

\subsection*{A.2. DRL for Routing Problems}
DRL-based methods for routing problems can be classified into two categories: constructive methods and improvement methods. Constructive methods focus on learning policies to construct solutions step-by-step. \cite{r_3} pioneered a sequence-to-sequence model, the Pointer Network (Ptr-Net), to solve the TSP in an autoregressive manner. Later, with the development of self-attention mechanism, \cite{r_4} introduced an attention-based model (AM), which achieved superior performance in various routing problems and became a milestone in this field. Since then, several AM-variants have been proposed, with most advanced neural solvers built on top of them \citep{r_38, r_36, r_59}. However, due to their limitations in the decision framework, problem representation, and training scheme, these approaches can only handle routing problems with simple structures and cannot be readily extended to more complex problems such as TPPs. Notably, while recent research has explored general models for different VRP variants \citep{r_60, r_61}, these efforts still remain confined to combinations of straightforward VRP configurations.

As another line of research, improvement methods leverage DRL to iteratively refine an initial solution through local search \citep{r_33, r_76}. It is worth mentioning that our approach is compatible with such improvement methods, and could potentially be integrated for further improvement in applications with larger computational budgets.
	
\subsection*{A.3. Meta-Learning}

Meta-learning, often referred to as "learning to learn", is a machine learning paradigm designed to enhance the adaptability and efficiency of learning algorithms across new tasks or data distributions by leveraging prior training experience \citep{r_77}. This idea of meta-learning  traces back decades ago, aiming to enable systems to adapt to new tasks with minimal data by leveraging prior knowledge, similar to how humans learn new skills based on past experiences \citep{r_63}. A significant milestone in meta-learning was the introduction of model-agnostic meta-learning (MAML) \citep{r_65}, which learns the initial parameters of neural networks such that it can be efficiently fine-tuned with a few gradient steps for fast adaptation to new tasks. Follow-up works have focused on improving the efficiency and robustness \citep{r_69}, and exploring other application domains such as RL \citep{r_68}. To date, the applications of meta-learning have expanded to complex tasks such as robotics \citep{r_66}, optimization \citep{r_64}, and multi-agent systems \citep{r_67}. For a comprehensive overview, readers are referred to the survey by \citep{r_77}.

\section{Ablation Study} \label{Ablation Study}
To get a better understanding of the contributions of components in our framework, we conduct a series of ablation studies. Following the organization of this paper, we study the contributions of each key component, including the bipartite graph representation, the policy network architecture, and the meta-learning strategy. The experimental results are presented in Table \ref{Ablation}, where we compare the performance of RL-E2E to exclude the effect of the post-optimization procedure. Besides, as the differences in runtime are very marginal, we focus on the objective value for comparative analysis.

\subsection*{B.1. Bipartite Graph Representation}
We evaluate the contributions of our bipartite graph representation by comparing it with the 1) complete graph and 2) \emph{k-NN} graph representations (entries ``Cplt.'' and ``\emph{k-NN}'' in Table \ref{Ablation}). The node features are defined as described in Section \ref{Bipartite Graph Representation for TPPs}, and the \emph{k-NN} graphs are adapted using attention masks in the encoder, with neighbor selection following \citep{r_5}. The results show that, while the complete and \emph{k-NN} graphs encode product supply information through node features, our bipartite graph representation facilitates the capture of high-level information for both markets and products through edges connecting the market nodes and product nodes, thus leading to better solutions. Notably, the \emph{k-NN} graph representation exhibited a significant performance drop, contrasting with its effectiveness in TSP. This is because the route construction for TPPs should take into account not only the spatial relations between markets but also their purchasing interdependencies. Consequently, \emph{k-NN} graphs, which emphasize spatial proximity, can potentially misguide the route construction and result in poor performance.

\begin{table}[!h]
\centering
{\color{black}
\setlength{\tabcolsep}{8pt}
\fontsize{8}{8}\selectfont
\renewcommand\arraystretch{1.2}
\centering
\caption{Ablation study on representation and network architecture} \label{Ablation}
\vspace{5pt}
\begin{threeparttable}
\begin{tabular}{ccc|ccccccccc}
\toprule
\multicolumn{3}{c}{Instance} & \multirow{2}{*}{\raisebox{-0.26cm}{\makecell{\;Full \\ \;Model}}} & \multicolumn{2}{c}{Repre.} & \multicolumn{3}{c}{Emb. Module} & \multicolumn{3}{c}{Dec. Context} \\
\cmidrule(lr){1-3} \cmidrule(lr){5-6} \cmidrule(lr){7-9} \cmidrule(lr){10-12}
$|M|$ & $|K|$ & $\lambda$ & & Cplt.\tnote{1} & \emph{k-NN}\tnote{2} & IEM\tnote{3} & ME\tnote{4} & LSTM\tnote{5} & GE\tnote{6} & DC\tnote{7} & RC\tnote{8} \\
\midrule
50  & 50  & /    & \textbf{1897} & 1941 & 1922 & 2985 & 1986 & 1911 & 1928 & 1933 & 1976 \\
50  & 100 & /    & \textbf{2542} & 2570 & 3521 & 3923 & 2712 & 2570 & 2583 & 2597 & 2630 \\
100 & 50  & /    & \textbf{1563} & 1599 & 1574 & 2897 & 1689 & 1579 & 1612 & 1667 & 1743 \\
100 & 100 & /    & \textbf{2111} & 2159 & 2834 & 4010 & 2296 & 2125 & 2165 & 2205 & 2338 \\
\midrule
50  & 50  & 0.99 & \textbf{2032} & 2084 & 2528 & 3085 & 2102 & 2052 & 2048 & 2048 & 2197 \\
50  & 100 & 0.99 & \textbf{2567} & 2865 & 3616 & 3865 & 2762 & 2581 & 2595 & 2611 & 2883 \\
100 & 50  & 0.99 & \textbf{1753} & 1799 & 2206 & 3067 & 1962 & 1778 & 1818 & 1888 & 2161 \\
100 & 100 & 0.99 & \textbf{2302} & 2434 & 3348 & 4111 & 2607 & 2338 & 2351 & 2450 & 2792 \\
\midrule
50  & 50  & 0.95 & \textbf{2645} & 2680 & 3501 & 3521 & 2857 & 2663 & 2714 & 2668 & 3065 \\
50  & 100 & 0.95 & \textbf{3457} & 3531 & 4930 & 4411 & 3645 & 3497 & 3488 & 3485 & 4171 \\
100 & 50  & 0.95 & \textbf{3027} & 3106 & 3536 & 4118 & 3451 & 3098 & 3138 & 3087 & 3711 \\
100 & 100 & 0.95 & \textbf{3993} & 4165 & 5629 & 5240 & 4461 & 4058 & 4073 & 4079 & 4884 \\
\midrule
50  & 50  & 0.9  & \textbf{3704} & 3763 & 4900 & 4552 & 3975 & 3744 & 3777 & 3769 & 4235 \\
50  & 100 & 0.9  & \textbf{4927} & 5047 & 6563 & 5737 & 5184 & 4985 & 4951 & 4967 & 5658 \\
100 & 50  & 0.9  & \textbf{4956} & 5110 & 5266 & 6016 & 5580 & 5047 & 5121 & 5036 & 6034 \\
100 & 100 & 0.9  & \textbf{6672} & 6869 & 7715 & 7823 & 7531 & 6869 & 6899 & 6838 & 8209 \\
\bottomrule
\end{tabular}

\begin{tablenotes}
\begin{minipage}[t]{0.4\textwidth}
\item[1] complete graph
\item[2] \emph{k-NN} graph
\item[3] w/o input embedding module
\item[4] w/o market encoder
\end{minipage}\hfill
\begin{minipage}[t]{0.4\textwidth}
\item[5] w/o LSTM module
\item[6] w/o global embedding
\item[7] w/o demand context
\item[8] w/o route context
\end{minipage}
\end{tablenotes}

\end{threeparttable}
}
\end{table}

\subsection*{B.2. Policy Network Architecture}
We study the contributions of the key components in our policy network, focusing on the embedding modules and the decoding context. Specifically, for the embedding modules, we evaluate the removal of 1) the input embedding module (IEM), 2) the market encoder (ME), and 3) the LSTM for current route embedding in the decoder. For the decoding context, we evaluate the removal of 1) the global embedding (GE), 2) the demand context (DC), and 3) the route context (RC). The results are presented in Table \ref{Ablation} under categories ``Emb. Module'' and ``Dec. Context''.

The results indicate that all three embedding modules contribute to overall performance, with the input embedding module being the most critical. Removing the input embedding module leads to a significant performance drop, as its removal separates the market locations from the product supply information in the input. This underscores the necessity of jointly considering spatial and supply information for TPPs. The market encoder becomes increasingly important with larger problem sizes, particularly for instances with more markets, owing to its capacity to extract deeper relational information between markets. Furthermore, while LSTM-based embeddings for the current partial route are not commonly employed in AM-based neural solvers \citep{r_4}, we find that this component enhances solution quality for TPPs, since all visited markets can influence subsequent decisions in solving the TPP.

Interestingly, the removal of global embedding has only a marginal impact on the solution quality. This suggests that the first one-to-many attention layer in the decoder can also serve as a mechanism for aggregating global information. Moreover, as discussed in Section \ref{State}, explicitly providing the demand context to the decoder also contributes to improved performance, especially when accurately meeting the demands plays a more important role in route construction. Finally, it is noteworthy that removing the route context, where our model becomes nearly non-autoregressive, still yields reasonable solutions. This suggests our policy network holds strong potential for non-autoregressive modeling---a parallelizable framework considered advantageous for scaling to extremely large instances \citep{r_70, r_71}. This finding opens up possibilities for further improvements in scalability in future work.

\subsection*{B.3. Meta-Learning Strategy}

The necessity and effectiveness of our meta-learning strategy for training stability have been discussed in Section \ref{Training Performance}. Here we focus on its impact on the zero-shot generalization ability, by comparing it against two commonly-used learning strategies: 1) transfer learning (TL) \citep{r_72}, where a policy network pre-trained on small instances is used as initialization for training on large instances, and 2) curriculum learning (CL) \citep{r_73}, which progressively increases the size of training instances, starting from smaller ones and gradually scaling to larger ones. The comparisons of zero-shot generalization performance are reported in Table \ref{meta in generalization}, demonstrating the advantage of our meta-learning strategy in generalizing to previously unseen, larger-sized instances.

\begin{table}[!h]
\centering
{\color{black}
\setlength{\tabcolsep}{8pt}
\fontsize{8}{8}\selectfont
\renewcommand\arraystretch{1.2}
\centering
\caption{Comparisons of zero-shot generalization} \label{meta in generalization}
\vspace{5pt}
\begin{threeparttable}
\begin{tabular}{ccc|cccccccc}
\toprule
\multicolumn{3}{c}{Instance} & \multicolumn{2}{c}{Heuristic} & \multicolumn{2}{c}{Meta} & \multicolumn{2}{c}{TL} & \multicolumn{2}{c}{CL} \\
\cmidrule(lr){1-3} \cmidrule(lr){4-5} \cmidrule(lr){6-7} \cmidrule(lr){8-9} \cmidrule(lr){10-11}
$|M|$ & $|K|$ & $\lambda$ & \raisebox{0.0cm}{\makecell{\fontsize{6}{6}\selectfont \; GSH \\[-0.5ex] \fontsize{6}{6}\selectfont +TRH}} & \raisebox{0.0cm}{\makecell{\fontsize{6}{6}\selectfont \; CAH \\[-0.5ex] \fontsize{6}{6}\selectfont +TRH}} & E2E\tnote{1} & TRH\tnote{2} & E2E & TRH & E2E & TRH \\
\midrule
150 & 150 & / & 2672  & 2460  & 2514  & \textbf{2380}  & 2745 & 2458 & 2648  & 2398  \\
200 & 200 & / & 2885  & 2454  & 2485  & \textbf{2262}  & 2850 & 2497 & 2806  & 2339  \\
300 & 300 & / & 3445  & \textbf{3074}  & 3206  & 3088  & 4367 & 3243 & 4043  & 3164  \\
\midrule
150 & 150 & 0.99 & 2918  & 3056  & 2576  & \textbf{2519}  & 2797 & 2604 & 2898  & 2694  \\
150 & 150 & 0.95 & \textbf{5642}  & 5884  & 6133  & 6014  & 6316 & 6174 & 6454  & 6224  \\
150 & 150 & 0.9  & \textbf{10008} & 10199 & 10707 & 10155 & 12061 & 10908 & 16027 & 12478 \\
\bottomrule
\end{tabular}

\begin{tablenotes}
\item[1] RL-E2E based on meta-learning strategy
\item[2] RL+TRH based on meta-learning strategy
\end{tablenotes}

\end{threeparttable}
}
\end{table}

\section{Extended Results} \label{Extended Results}

At last, to further highlight the effectiveness and compatibility of our approach, we present extended results that include comparisons with state-of-the-art neural VRP solvers, as well as a validation of the cross-problem compatibility of our approach and network architecture with other VRP variants.

\subsection*{C.1. Comparisons with Other Neural Solvers}

We compare our approach with several neural solvers, including the widely-recognized Ptr-Net \citep{r_3} and AM \citep{r_4}, as well as the state-of-the-art unified neural VRP solver, MVMoE \citep{r_61}. Since the baseline neural VRP solvers are not readily applicable to TPPs, we adapt their node features as described in Section \ref{Bipartite Graph Representation for TPPs}, and employ transfer learning to mitigate potential training collapse on large-sized instances. The results are presented in Table \ref{other neural solvers}. Our approach consistently outperforms the above neural solvers across all instance sets.

Beyond superior solution quality, our approach offers a significant advantage in scalability. The baseline neural solvers, limited by their problem representations, require separate models and training procedures for TPP instances with varying numbers of products. Even MVMoE, despite introducing a unified policy network with feature padding to multiple VRP variants, remains restricted by a fixed-dimensional feature representation, limiting its flexibility across TPP instances of varying sizes. This limitation not only hinders them from learning generalized knowledge but also signiﬁcantly increases the training and storage overhead in practice. In contrast, our bipartite graph representation enables a size-agnostic model, allowing us to train a single, generalizable policy network through meta-learning that effectively adapts across varied instance sizes and distributions, significantly enhancing its flexibility and efficiency for real-world applications.

\begin{table}[!h]
\centering
{\color{black}
\setlength{\tabcolsep}{8pt}
\fontsize{8}{8}\selectfont
\renewcommand\arraystretch{1.2}
\centering
\caption{Comparisons with other neural solvers}  \label{other neural solvers}
\vspace{5pt}
\begin{tabular}{ccc|cccccccc}
\toprule
\multicolumn{3}{c}{Instance} & \multicolumn{2}{c}{Ptr-Net} & \multicolumn{2}{c}{AM} & \multicolumn{2}{c}{MVMoE} & \multicolumn{2}{c}{Ours}\\
\cmidrule(lr){1-3} \cmidrule(lr){4-5} \cmidrule(lr){6-7} \cmidrule(lr){8-9} \cmidrule(lr){10-11}
$|M|$ & $|K|$ & $\lambda$ & E2E & TRH & E2E & TRH & E2E & TRH & E2E & TRH \\
\midrule
50  & 50  & /    & 2427 & 2056 & 2009 & 1906 & 2210 & 1976 & 1897 & \textbf{1857} \\
50  & 100 & /    & 3357 & 2677 & 2856 & 2525 & 2813 & 2507 & 2542 & \textbf{2446} \\
100 & 50  & /    & 2284 & 2027 & 1643 & 1556 & 1938 & 1727 & 1563 & \textbf{1524} \\
100 & 100 & /    & 3114 & 2394 & 2404 & 2148 & 2476 & 2186 & 2111 & \textbf{2044} \\
\midrule
50  & 50  & 0.99 & 2997 & 2430 & 2133 & 2005 & 2272 & 2049 & 2032 & \textbf{1954} \\
50  & 100 & 0.99 & 3891 & 2821 & 3093 & 2650 & 2817 & 2546 & 2567 & \textbf{2466} \\
100 & 50  & 0.99 & 2284 & 2027 & 1863 & 1777 & 2123 & 1927 & 1753 & \textbf{1711} \\
100 & 100 & 0.99 & 3598 & 2867 & 2720 & 2465 & 2665 & 2421 & 2302 & \textbf{2235} \\
\midrule
50  & 50  & 0.95 & 3281 & 2925 & 2843 & 2678 & 2922 & 2723 & 2645 & \textbf{2594} \\
50  & 100 & 0.95 & 4255 & 3788 & 3857 & 3529 & 3714 & 3469 & 3457 & \textbf{3368} \\
100 & 50  & 0.95 & 3744 & 3451 & 3280 & 3137 & 3490 & 3284 & 3027 & \textbf{2980} \\
100 & 100 & 0.95 & 4996 & 4523 & 4400 & 4162 & 4465 & 4218 & 3993 & \textbf{3920} \\
\midrule
50  & 50  & 0.9  & 4391 & 4014 & 3938 & 3758 & 4024 & 3819 & 3704 & \textbf{3644} \\
50  & 100 & 0.9  & 5676 & 5257 & 5305 & 5066 & 5289 & 5025 & 4927 & \textbf{4855} \\
100 & 50  & 0.9  & 5864 & 5488 & 5282 & 5134 & 5623 & 5313 & 4956 & \textbf{4900} \\
100 & 100 & 0.9  & 7664 & 7246 & 7223 & 6946 & 7426 & 7018 & 6672 & \textbf{6660} \\
\bottomrule
\end{tabular}
}
\end{table}

\subsection*{C.2. Compatibility with Other VRP Variants}
It is worth mentioning that many VRP variants can be formulated as special cases of the TPP. For example, the TSP is equivalent to the UTPP where each market provides a unique product with zero price. Similarly, the prizes and penalties in routing problems such as the orienteering problem (OP) \citep{r_74} and the prize collecting TSP (PCTSP) \citep{r_75} can also be represented as products with associated purchasing costs in TPPs. Therefore, our approach is inherently compatible with these VRP variants. To validate this, we evaluate our approach on the TSP, OP, and PCTSP, following the experimental settings in \citep{r_4}. The results are summarized in Table \ref{other VRP variants}. It is demonstrated that our approach outperforms AM on most problems, mainly due to our improvements in the embedding module, decoding context, and training scheme. Moreover, our approach could potentially be further enhanced by incorporating training and inference techniques tailored to these specific VRP variants.

\begin{table}[!b]
\centering
{\color{black}
\setlength{\tabcolsep}{8pt}
\fontsize{8}{8}\selectfont
\renewcommand\arraystretch{1.2}
\centering
\caption{Compatibility with other VRP variants}  \label{other VRP variants}
\vspace{5pt}
\begin{threeparttable}
\begin{tabular}{c|ccc|ccc|ccc}
\toprule
\multirow{2}{*}{Method} & \multicolumn{3}{c}{TSP ($\downarrow$)\tnote{*}} & \multicolumn{3}{c}{OP ($\uparrow$)\tnote{*}} & \multicolumn{3}{c}{PCTSP ($\downarrow$)\tnote{*}} \\
& $n$=20 & $n$=50 & $n$=100 & $n$=20 & $n$=50 & $n$=100 & $n$=20 & $n$=50 & $n$=100 \\
\midrule
AM & 3.85 & 5.80 & 8.12 & 5.19 & 15.64 & \textbf{31.62} & 3.18 & 4.60 & 6.25 \\
Ours & \textbf{3.84} & \textbf{5.73} & \textbf{7.97} & \textbf{6.54} & \textbf{16.20} & 31.49 & \textbf{3.04} & \textbf{4.48} & \textbf{6.11} \\
\bottomrule
\end{tabular}

\begin{tablenotes}
\item[*] $\downarrow$ means a minimization problem, and $\uparrow$ means a maximization problem
\end{tablenotes}

\end{threeparttable}
}
\end{table}

\end{document}